\documentclass{article}
\usepackage{amsmath,amsfonts,latexsym, amsrefs,amssymb}

\usepackage{graphicx}

\DeclareMathOperator{\Prob}{\mathbb{P}}   

\newcommand{\1}{\mathds{1}}
\usepackage{enumerate}

\usepackage{wrapfig}

\usepackage{color}

\numberwithin{equation}{section}

\newcommand{\rd}{{\rm d}}
\newcommand{\rdA}{{\rm d A}}

\newcommand{\by}{{\bf{y}}}

\newcommand{\mU}{{\mathbb{U}}}

\newcommand{\mT}{{\mathbb{T}}}

\newcommand{\be}{\begin{equation}}
\newcommand{\ee}{\end{equation}}

\newcommand{\hp}[1]{with $ #1 $-high probability}

\newcommand{\e}{{\varepsilon}}

\newcommand{\T}{\mathbb T}

\newcommand{\non}{\nonumber}

\usepackage{amsmath} 
\usepackage{amssymb}
\usepackage{amsthm}
\usepackage{amsxtra}
\usepackage{amscd}
\usepackage{bbm}
\usepackage{mathrsfs}
\usepackage{bm}

\newcommand{\bb}{\mathbb} 

\renewcommand{\cal}{\mathcal}

\newcommand{\wt}{\widetilde}
\newcommand{\mG}{\mathcal G}

\newcommand{\ii}{\mathrm{i}} 

\newcommand{\deq}{\mathrel{\mathop:}=}

\renewcommand{\epsilon}{\varepsilon}
\renewcommand{\leq}{\leqslant}
\renewcommand{\geq}{\geqslant}

\renewcommand{\le}{\leq}
\renewcommand{\ge}{\geq}

\renewcommand{\P}{\mathbb{P}}
\newcommand{\E}{\mathbb{E}}
\newcommand{\R}{\mathbb{R}}
\newcommand{\C}{\mathbb{C}}
\newcommand{\N}{\mathbb{N}}
\newcommand{\Z}{\mathbb{Z}}

\newcommand{\htP}{\mathscr{P}}

\newcommand{\qB}[1]{\Bigl[{#1}\Bigr]}

\newcommand{\absb}[1]{\bigl\lvert #1 \bigr\rvert}

\DeclareMathOperator{\tr}{Tr}

\DeclareMathOperator{\supp}{supp}

\DeclareMathOperator{\re}{Re}
\DeclareMathOperator{\im}{Im}

\DeclareMathOperator{\OO}{O}
\DeclareMathOperator{\oo}{o}

\theoremstyle{plain} 
\newtheorem{theorem}{Theorem}[section]
\newtheorem*{theorem*}{Theorem}
\newtheorem{lemma}[theorem]{Lemma}
\newtheorem*{lemma*}{Lemma}
\newtheorem{corollary}[theorem]{Corollary}
\newtheorem*{corollary*}{Corollary}

\newtheorem*{proposition*}{Proposition}

\newtheorem{definition}[theorem]{Definition}
\newtheorem*{definition*}{Definition}

\newtheorem*{eximple*}{Eximple}

\newtheorem*{remark*}{Remark}
\newtheorem*{remarks*}{Remarks}

\makeatletter
\renewcommand{\section}{\@startsection
{section}
{1}
{0mm}
{-2\baselineskip}
{1\baselineskip}
{\normalfont\large\scshape\centering}} 
\makeatother

\makeatletter
\renewcommand{\subsection}{\@startsection
{subsection}
{2}
{0mm}
{-\baselineskip}
{0\baselineskip}
{\normalfont\bf} } 
\makeatother

\usepackage[T1]{fontenc}

\usepackage{dsfont}
\usepackage{stmaryrd}

\newcommand{\cor}{\color{red}}

\newcommand{\nc}{\normalcolor}

\begin{document}

\title{{\sc\Large The local circular law III: general case}\vspace{1cm}}

\date{}

\author{\vspace{0.5cm}
{\sc Jun Yin}$
$\thanks{Partially supported
by NSF grant DMS-1001655 and DMS- 1207961}
 \\\\
\normalsize Department of Mathematics, University of Wisconsin-Madison \\
\normalsize Madison, WI 53706-1388, USA \ \normalsize jyin@math.wisc.edu $ $\vspace{1cm}}

\maketitle

\begin{abstract}
In the first part \cite{BouYauYin2012Bulk } of this   article series,  Bourgade, Yau and the author of this paper proved a local version of the circular law up to the finest scale $N^{-1/2+ \e}$
for  non-Hermitian random matrices at any point $z \in \C$ with $||z| - 1| > c $ for any constant $c>0$ independent of
the size of the matrix. In the second part  \cite{BouYauYin2012Edge}, they extended this result to include the edge case  $ |z|-1=\oo(1)$, under the main assumption that the  third moments of the matrix elements vanish. (Without the vanishing third moment assumption, they proved that the circular law is valid near the spectral edge $ |z|-1=\oo(1)$
up to scale $N^{-1/4+ \e}$.)  In this paper, we will remove this assumption, i.e. we prove a local version of the circular law up to the finest scale $N^{-1/2+ \e}$
for  non-Hermitian random matrices at any point $z \in \C$. 
 \end{abstract}

\vspace{1.5cm}


\medskip

\medskip


\medskip

\newpage

   \section{Introduction and Main result}
 
 The circular law in random matrix theory describes the macroscopic limiting spectral measure of normalized non-Hermitian matrices with
independent entries. Its origin goes beck to the work of Ginibre \cite{Gin1965},
who found the joint density of the eigenvalues of such Gaussian matrices. More precisely, for an
$N\times N$ matrix with independent  entries $\frac{1}{\sqrt {N}}z_{ij}$ such that  $z_{ij}$ is
identically distributed according to the measure $\mu_g=\frac{1}{\pi}e^{-|z|^2}\rdA(z)$
($\rdA$ denotes the Lebesgue measure on $\mathbb{C}$),
its eigenvalues $\mu_1,\dots,\mu_N$ have a probability density
proportional to
\begin{equation}\label{eqn:Ginibre}
\prod_{i<j}|\mu_i-\mu_j|^2e^{-N\sum_{k}|\mu_k|^2}
\end{equation}
with respect to the Lebesgue measure on $\C^{N}$.
These random spectral measures define a determinantal point process with the
explicit kernel (see \cite{Gin1965}) 
\begin{equation}\label{eqn:GinibreKernel}
K_N(z_1,z_2)=\frac{N}{\pi}e^{-\frac{N}{2}(|z_1|^2+|z_2|^2)}\sum_{\ell=0}^{N-1}\frac{(N z_1\overline{z_2})^\ell}{\ell!}
\end{equation}
with respect to the Lebesgue measure on $\C$. This integrability property allowed Ginibre to derive the  circular law for the eigenvalues,
i.e.,
$\frac{1}{N}\rho_1^{(N)}$ converges to the uniform measure on the unit circle,
\begin{equation}\label{eqn:circular}
\frac{1}{\pi}\1_{|z|<1}\rdA(z).
\end{equation}
This limiting law also holds for real Gaussian entries \cite{Ede1997},
for which a more detailed analysis was performed in \cites{ForNag2007,Sin2007,BorSin2009}.

For non-Gaussian entries,
Girko \cite{Gir1984} argued that the macroscopic limiting spectrum is still given by
(\ref{eqn:circular}).  His main insight is commonly known as  the {\it Hermitization technique}, which
converts the convergence of complex empirical measures into the convergence of logarithmic transforms of
a family of Hermitian matrices.  If we denote the original non-Hermitian matrix by $X$ and the eigenvalues of
$X$ by $\mu_j$,  then for any ${C}^2$ function $F$ we have the identity
\be\label{id0}
\frac 1 N \sum_{j=1}^N F (\mu_j) = \frac1{4\pi N} \int \Delta F(z)   \tr   \log (X^* - z^* ) (X-z)    \rdA(z).
\ee
Due to the logarithmic singularity at $0$, it is clear that the small eigenvalues of the Hermitian matrix $(X^* - z^* ) (X-z)  $
play a special role. A key question is to
estimate  the small eigenvalues of $(X^* - z^* ) (X-z)$, or in other words, the small singular values
of $ (X-z)$.
 This problem
was not treated  in \cite{Gir1984},
but the  gap was remedied in a series of  papers.
First Bai \cite{Bai1997} was able to treat the  logarithmic  singularity assuming  bounded density and bounded high moments
for the entries of the matrix (see also \cite{BaiSil2006}).
Lower bounds on the smallest singular values were given in Rudelson, Vershynin \cites{Rud2008,RudVer2008}, and subsequently
Tao, Vu \cite{TaoVu2008}, Pan, Zhou \cite{PanZho2010} and G\"otze, Tikhomirov \cite{GotTik2010} weakened the
moments and smoothness assumptions for the circular law, till the optimal $\mbox{L}^2$ assumption,
under which the circular law was proved in \cite{TaoVuKri2010}. On the other hand, Wood \cite{Wood2010} showed that the  circular law also holds
 for sparse
random $n$ by $n$ matrices where each entry is nonzero with probability $ n^{\alpha-1}$ where $0<\alpha\leq 1$. 

\bigskip

In the first part of this article \cite{BouYauYin2012Bulk}, Bourgade, Yau and the author of this paper proved a local version of the circular law, up to the optimal scale
$N^{-1/2 + \e}$, in the bulk of the spectrum.  In the second part  \cite{BouYauYin2012Edge}, they extended this result to include the edge case, under the assumption that the  third moments of the matrix elements vanish. (Without the vanishing third moment assumption, they also  proved that the circular law is valid near the spectral edge $ |z|-1=\oo(1)$
up to scale $N^{-1/4+ \e}$.) This vanishing third moment condition is  also the main assumption  in  Tao and Vu's  work on local circular law \cite{TaoVu2012}.  In the current paper, we will remove this assumption, i.e. we prove a local version of the circular law up to the finest scale $N^{-1/2+ \e}$
for  non-Hermitian random matrices at any point $z \in \C$. 

More precisely,
we considered   an $N \times N$  matrix $X$ with independent real\footnote{For the sake of notational simplicity we do not consider complex entries in this paper, but the statements and proofs are similar.} centered entries with variance $ N^{-1}$.
Let $\mu_j$, $j\in \llbracket 1, N\rrbracket$ denote the eigenvalues of $X$. To state the local circular law, we first define the
notion of \emph{stochastic domination}. \begin{definition}\label{stodom} Let $W=W^{(N)}$ be a  family of random variables and $\Psi=\Psi^{(N)}$  be a family of    deterministic  parameters.
We say that $W$ is  \emph{stochastically dominated} by  $\Psi$
if for any   $  \sigma> 0$ and $D > 0$ we have
\begin{equation}\label{stodo}
 \P \qB{\absb{W } > N^\sigma \Psi  } \;\leq\; N^{-D}
\end{equation}\nc
for sufficiently large $N$. 
We denote this stochastic domination property by
\begin{equation*}
W \;\prec\; \Psi\,,\quad or \quad W =\OO_\prec (\Psi).
\end{equation*}
Furthermore,  Let  $U^{(N)}$ be  a possibly N-dependent parameter set.  
We say  $W(u)$ is  \emph{stochastically dominated} by  $\Psi(u)$ uniformly in $u\in U^{(N)}$,  
if for any   $  \sigma> 0$ and $D > 0$ we have
\begin{equation}\label{stodo2}
\sup_{u\in U^{(N)}} \P \qB{\absb{W (u)} > N^\sigma \Psi (u) } \;\leq\; N^{-D}
\end{equation}\nc
for  uniformly   sufficiently large $N$ (may depends on $\sigma$ and $D$). 
\end{definition}
 Note: In the most cases of this paper, the $U^{(N)}$ is chosen as the product of the  index sets $1\leq i, j\leq N$ and some  compact set in $\C^2$.

\bigskip

In this paper, as in \cite{BouYauYin2012Bulk}, \cite{BouYauYin2012Edge} and \cite{TaoVu2012}, we assume that  the probability distributions of the   matrix elements
satisfy the following uniform subexponential decay property:
\be\label{subexp}
\sup_{(i,j)\in\llbracket 1,N\rrbracket^2}\Prob\left(|\sqrt{N}X_{i,j}|>\lambda\right)\leq \vartheta^{-1} e^{-\lambda^\vartheta }
\ee
for some constant $\vartheta >0$ independent of $N$.
This condition can of course be weakened to an hypothesis of boundedness on sufficiently high moments, but the error estimates
in the following Theorem would be weakened as well.

Note: most   constants appearing in this work may depend on $\vartheta$, but we will not emphasize this dependence in the proof.  

Let $f:\mathbb{C}\to\mathbb{R}$ be a fixed smooth compactly supported function,
and $f_{z_0}(\mu)=N^{2s}f(N^s(\mu-z_0))$, where    $z_0$ depends on $N$, and $s$ is a fixed scaling parameter in $[0,1/2]$. Let $D$ denote the unit disk.
Theorem 2.2 of  \cite{BouYauYin2012Bulk} and Theorem 1.2 of \cite{BouYauYin2012Edge} assert that  the following estimate holds: (Note: Here $\|f_{z_0}\|_1=O(1)$)
\begin{equation}\label{yjgq}
\left( N^{-1} \sum_{j}f_{z_0} (\mu_j)-\frac1\pi \int_D f_{z_0}(z)    \, \rdA(z)  \right)\prec \ N^{-1+2s },  \quad s\in (0,1/2],
\end{equation}
 \cite{BouYauYin2012Bulk}  if   $||z_0|-1|>c$ for some $c>0$ independent of $N$ or  \cite{BouYauYin2012Edge} if the third moments of matrix entries vanish. 
This implies that the circular law holds after zooming up to scale $N^{-1/2+\e}$ ($\e>0$) under these conditions.
In particular, there are neither clusters of eigenvalues nor holes in the spectrum at such scales. We note that in \cite{BouYauYin2012Bulk} and 
\cite{BouYauYin2012Edge}, the scaling parameter was denoted as $a$, but the letter $a$ will be used as a fixed index in this work. 

We aim at understanding the circular law for any $z_0\in \C$ without the vanishing third moment assumption.   The following theorem is our main result.

\begin{theorem}\label{z1} {\bf Local circular law:}
Let  $X$ be an $N\times N$ matrix with independent centered
entries of  variances $1/N$.   Suppose that
the distributions of the matrix elements   satisfy the subexponential decay property (\ref{subexp}).
Let $f_{z_0}$ be defined as above \eqref{yjgq} and $D$ denote the unit disk.
Then  for any $s \in (0,1/2]$ and any $z_0\in \C$, we have
\be\label{yjgq2}
 \left( N^{-1} \sum_{j}f_{z_0} (\mu_j)-\frac1\pi \int_D f_{z_0}(z)\rdA(z)     \right)\prec \ N^{-1+2s }.
\ee
Notice that the main new assertion of \eqref{yjgq2}  is for  the case: $|z_0|-1=\oo(1)$ and the third moments not vanishing, since the other cases were proved in  
\cite{BouYauYin2012Bulk} and \cite{BouYauYin2012Edge},
stated in \eqref{yjgq}.
\end{theorem}

Remark: Shortly after the preprint \cite{BouYauYin2012Bulk} appeared, a version of
local circular law (both in the bulk and near the edge) was proved by  Tao and Vu \cite{TaoVu2012} under the assumption
that  the first three moments of the matrix entries match a Gaussian distribution, i.e., the third moment vanish.

In the next section we will introduce our main strategy and improvements. 

 \section{Proof of Theorem \ref{z1}}

{\it Proof of Thm. \ref{z1}: } The bulk case of  Thm. \ref{z1}  was proved in Theorem 2.2 of  \cite{BouYauYin2012Bulk}. Furthermore, it is easy to see that the results in Thm. \ref{z1} for $s=1/2$  follow from the results in  for $s<1/2$.  Hence in this proof, we can assume that 
 $$
 ||z_0|-1|=o(1), \quad s\in (0,1/2)
 $$
In the  edge case,  our Thm. \ref{z1} was proved in the Thm 1.2 of \cite{BouYauYin2012Edge} with the vanishing third moment assumption.  Hence the goal of this paper is to improve the proof of Thm. 1.2 of \cite{BouYauYin2012Edge}. One can easily check that in the proof of  Thm. 1.2 of \cite{BouYauYin2012Edge}, the condition $\E X^3_{ij}=0$  was only used in the Lemma 2.13 of \cite{BouYauYin2012Edge}.  Therefore,  we only 
 need to prove a stronger version of Lemma 2.13  in  \cite{BouYauYin2012Edge}, i.e., the one without vanishing third moment condition. More precisely, 
 it only remains to prove  the following lemma \ref{lyy}. (Here we use the same notations as   in \cite{BouYauYin2012Edge}, except for the scaling parameter)
 
 \qed

  Before stating   lemma \ref{lyy}, i.e., the stronger version of Theorem 1.2 of \cite{BouYauYin2012Edge}, we introduce some definitions and notations.  
  First, we introduce  the notation
$$
 Y:=Y_z:=X-z I  \quad
$$
 where $I$ is the identity operator. In the following, we use the notation $A\sim B$ when $c B \le  |A|\leq c^{-1}B$,
where $c>0$ is independent of $N$. For any  matrix $M$, we denote $M^T$ as the transpose of $M$ and 
$M^*$ as the Hermitian conjugate.  Usually we choose $z -z_0\sim N^{-s}$,  hence we define the scaled parameter $\xi$:
 $$
 z=z_0+N^{-s}\xi,\; i.e., \; \xi:= N^s(z-z_0)
 $$
Define the Green function
 of $Y^*_z Y_z$  and its trace by, where $w\in \C$ and $\im w>0$,
\be\label{defG}
 G(w):=G(w,z)=(Y^*_z Y_z-w)^{-1},\quad   m(w):=m(w,z)=\frac{1}{N}\tr G(w,z)=\frac{1}{N}\sum_{j=1}^N\frac{1}{ \lambda_j(z) - w}.
\ee
 Let  $ m_{\rm c}:=m_{\rm c}(w,z)$  be  the unique
solution of
 \be\label{defmc1}
 m_{\rm c}^{-1}=-w(1+m_{\rm c})+|z|^2(1+m_{\rm c})^{-1}
\ee
with positive  imaginary  part. As proved in \cite{BouYauYin2012Bulk} and \cite{BouYauYin2012Edge},  for some regions of $(w, z)$ with high  probability, 
$m(w,z)$ converges to  $m_{\rm c}(w,z)$ pointwise,
as $N\to \infty$.  Let $\rho_c$ be the measure whose Stieltjes transform is $m_{\rm c} $.  This measure is compactly supported and $\supp \rho_c=[\max\{0, \lambda_-\}, \lambda_+]$, where 
 \be\label{deflampm}
 \lambda_\pm :=\lambda_{\pm }(z):=\frac{( \alpha\pm3)^3}{8(\alpha\pm1)} ,\quad \alpha:=\sqrt{1+8|z|^2}.
\ee
Note  that $\lambda_-$ has the same sign as $|z|-1$.
It is well-known that $ \rho_{\rm c} (x,z)$   can be obtained from its Stieltjes transform $m_{\rm c}(x+\ii\eta,z)$ via
$$
 \rho_{\rm c} (x,z)= \frac1\pi\im  \lim_{\eta\to 0^+}m_{\rm c}(x+\ii\eta,z)
=\frac1\pi\1_{x\in[\max\{0,\lambda_-\},\lambda_+]} \im  \lim_{\eta\to 0^+}m_{\rm c}(x+\ii\eta,z).
$$
(Some  basic properties of $m_c$ and $\rho_c$ were discussed in  section 2.2 of \cite{BouYauYin2012Edge})
 \begin{definition}\label{defphiY}  { $\phi$, $\chi$, $I$ and $Z_{X,{\rm c}}^{(f)}$}

Let $h(x)$ be a smooth  increasing  function supported on $[1,+\infty]$  with  $h(x)=1$
for $x\geq 2 $ and $h(x)=0$  for $x\leq 1$.
For any $\e>0$, define $\phi$ on $\R_+$ by (note: $\lambda_+$ depends on $z$)
\be\label{defvarphi}
 \phi(x):=\phi_{\e,z}(x):= h(N^{2-2\e} x)\, (\log x)\,\left(1- h\left(\frac x{2\lambda_+ }\right)\right).
 \ee
 Let $ \chi$  be a smooth cutoff function  supported in $ [-1,1]$ with  bounded derivatives and $ \chi(y) = 1$ for $|y| \le 1/2$. 
 Recall $\rdA$ denotes the Lebesgue measure on $\mathbb{C}$, for any fixed function $g$ defined on $\C$, we define:
 $$Z_{X,{\rm c}}^{(g)}:=Z_{X,{\rm c} }^{(g)}(z_0, \e, s) :=N
\int  \Delta g(\xi)  \int_{I}  \chi(\eta)\phi '(E)
\re (m(w)-m_{\rm c}(w)) \rd E\rd \eta\rdA(\xi), \quad w=E+i\eta,\quad  z=z_0+N^{-s}\xi $$
and \be\label{defa1a2}
I:=I_\e :=\left\{ w\in \C: N^{-1+\e}\sqrt E\leq \eta, \; E\ge N^{-2+2\e}, \;\; |w|\leq \e,\; w=E+i\eta\right\}.
\ee
\end{definition}
Note: the condition $E\ge N^{-2+2\e}$ was not in the definition of the $I$  used in \cite{BouYauYin2012Edge}, but clearly this condition is implied by $\phi'(E)\neq 0$, i.e., 
our new $I$ does not change the 
value of $Z_{X,{\rm c}}^{(g)}$. One can also easily check:
\be\label{wzzln}
w\in I_\e\implies |w|^{1/2}\leq 2N^{1-\e}\eta
\ee

With these notations and definitions, we claim the following  main lemma. It is a stronger version of Lemma 2.13  in  \cite{BouYauYin2012Edge}, i.e., the one without vanishing third moment condition.
\begin{lemma}\label{lyy} Under the assumptions of Theorem \ref{z1}, there exists a constant $C>0$ such that for any small
enough $\e>0$(independent of $N$), if $||z_0|-1|\leq \e$, $s\in (0, 1/2)$, then 
$$
Z_{X,{\rm c}}^{(f)}\prec N^{C\e}c_f,
$$
where $ c_f$ is a constant depending only on the function $f$.
\end{lemma}
As mentioned above,  in the proof of  Thm. 1.2 of \cite{BouYauYin2012Edge}, the vanishing third moment condition  was only used in the Lemma 2.13 of \cite{BouYauYin2012Edge}. Therefore with the  improved  Lemma \eqref{lyy}, one can obtain our main result theorem \ref{z1} as in \cite{BouYauYin2012Edge}.  

\qed

 In the next step,   the lemma \ref{lyy} will be reduced  to lemma \ref{main lemma}.
  
 \medskip

We note that the bounds  proved in \cite{BouYauYin2012Edge} for $G_{ij}$'s are not strong enough for our purpose  in this paper. Unfortunately we noticed that it seems impossible to  improve these bounds in general cases. 
 On the other hand,  we found that though the  behaviors $G$'s and $\mG$'s are unstable in the region $|m|\leq (N\eta)^{-1}$, they are very stable in the region $|m|\gg (N\eta)^{-1}$ and many stronger bounds can be derived in this region.  Therefore, in the following proof, we separate the $Z_{X,c}$ into two parts: the one comes for the  region $|m|\leq (N\eta)^{-1}$ and  the one comes for the  region $|m|\gg (N\eta)^{-1}$. The first part can be easily bounded, since the $m$ is small, so as its contribution  to $Z_{X,c}$.  For the second part, we will apply Green's function comparison method (which was first introduced in \cite{ErdYauYin2010PTRF} for generalized Wigner matrix) and our new stronger bounds  in the  region $|m|\gg (N\eta)^{-1}$.  
 
  On the other hand, the old Green's function comparison method was not enough for our purpose, which is also the reason that in \cite{BouYauYin2012Edge}, the authors needed the extra assumption on the third moment of the matrix entries.  In this work, we will introduce an improved Green's function comparison method, which provides an extra $N^{-1/2} $ factor than the previous method. This idea was motivated from the work in \cite{BEKYY2012}.
 \begin{definition}\label{DefAX} {$t_X$ and $A_X^{(f)}$}
 
 For $N\times N$ matrix $X$, we define
\be\nonumber
t_X:= t_X(  \e, w, z):=N^{-\e}N\eta\re m 
\ee
i.e.,
\be\nonumber
t_X:= N^{-\e}\eta\re\tr \left((X^*-z^*)(X-z)-w\right)^{-1}, \quad \eta=\im w
\ee
   {  Now we extend the function $h$ defined in Def. \ref{defphiY} to the whole real lane, i.e., $h(x)=h(-x)$, but still use the same notation $h(x)$.} With these notations, we define:
\be\label{defAXf}
A_{X }^{(f)}: =A_{X }^{(f)}(z_0, \e, s)=N
\int  \Delta f(\xi)  \int_{I}  \chi(\eta)\phi '(E)
\Big(h(t_X)\re  m-\re m_{\rm c}\Big) \rd E\rd \eta\rdA(\xi), \quad  
\ee
where $z=z_0+N^{-s}\xi$, $w=E+i\eta$,  $\phi=\phi_{\e,z}$ and $t_X= t_X(  \e, w, z)$. 
\end{definition}

 Note the only difference between $A_X^{(f)}$ and $Z_{X,{\rm c}}^{(f)}$ is the $h(t_X)$ in front of $\re m$. Then the difference of $A_X^{(f)}$ and $Z_{X,{\rm c}}^{(f)}$ only comes from  the region $h(t_X)\neq 1$, i.e, $|\re m| \leq 2N^{\e} (N\eta)^{-1}$. Therefore, by the definitions of $\phi$ we have 
 \be\label{Tgg}
 |A_{X }^{(f)}-Z_{X,{\rm c}}^{(f)}|\leq \int | \Delta f(\xi)|  \int_{I}  \chi(\eta)|\phi '(E)|
\left(2N^{\e} (N\eta)^{-1}\right) \rd E\rd \eta\rdA(\xi)\leq N^{C\e}c_f
 \ee
where we used $|(1-h(t_X))\re m|\leq 2 N^\e(N\eta)^{-1}$. 

\bigskip

{\it  Proof of Lemma \ref{lyy}: }
 With \eqref{Tgg}, it only remains to prove the following lemma. 
\begin{lemma}\label{main lemma} Under the assumptions of Theorem \ref{z1},   there exists a constant $C>0$ such that for any small
enough $\e>0$(independent of $N$), if $||z_0|-1|\leq \e$ and $  s\in(0, 1/2)$,  then  
\be\nonumber
A_{X }^{(f)} \prec N^{C\e}c_f
\ee
where $ c_f$ is a constant depending only on the function $f$.
 \end{lemma}
\qed

In the next subsection, we will introduce the basic idea of proving Lemma \ref{main lemma}. The rigorous proof will start from section 3. 
 
 \subsection{Basic strategy of proving  Lemma \ref{main lemma}: }\label{subsection: bs}

Before we give the complete proof of this lemma, we introduce the basic idea and main improvement in the remainder of this section.  Lemma \ref{lyy}  was proved in  \cite{BouYauYin2012Edge} under the vanishing third moment condition. With \eqref{Tgg}, that result implies  that if $X_{ij}$'s are Gaussian variables, for all $1\leq i,j\leq N$,  then  for any fixed $p\in 2\N$, 
\be\label{Bn1}
\E |A_{X,{\rm c}}^{(f)}|^p\prec N^{C\e p}, \quad\quad X_{ij}\sim \cal N(0, 1/N)
\ee
  As one can see that $A_{X,c}^{(f)}$ is  basically a linear functional of $m(w,z)$.  Hence  as in \cite{BouYauYin2012Edge}, we will apply  the Green function comparison method to show  that for    sufficiently large $N$, 
\be\label{Bn2}
 \E |A_{X}^{(f)}|^p  \leq C\,\E |A_{X'}^{(f)}|^p+N^{C\e p},  
\ee
for any two   different ensembles $X$ and $X'$ whose  matrix elements satisfy the condition of  Theorem \ref{z1}. To complete the proof for lemma \ref{main lemma}, we will choose $X'$ to be    the Ginibre ensemble, whose matrix elements are Gaussian variables. The $X$ will be the general ensembles in lemma \ref{main lemma}.   Combining \eqref{Bn1} and \eqref{Bn2}, with Markov inequality, one immediately obtains Lemma \ref{main lemma}. 

 In applying the Green function comparison method, we  estimate the expectation value of the functionals  of $Y $, $G=(Y^*  Y -w)^{-1}$  and $\mG=(YY^*   -w)^{-1}$, i.e., $\E F(Y, G, \mG)$.   In \cite{BouYauYin2012Edge} and most previous applications of  Green function comparison method, one can only  bound the expectation value of these functionals with their stochastically dominations.  For example, in  \cite{BouYauYin2012Edge}, for $i\neq j$ and $|w|^{1/2}\ll (N\eta)$, one has
$$
|(Y G)_{ij}|\prec 1
$$
With this stochastically domain, the authors   in \cite{BouYauYin2012Edge}  obtained that  $|\E (Y G)_{ij}|\leq N^{\sigma} $ for any $\sigma>0$. In the present paper, under the condition $|\re m|\gg (\eta N)^{-1}$, i.e., $h(t_X)>0$, we will first show an improved bound: for $i\neq j$ and $|w|^{1/2}\ll (N\eta)$
 $$
|h(t_X)(Y G)_{ij}|\prec  \sqrt{\frac{|w|^{1/2}}{N\eta} }
$$
Then using a new idea on Green's function comparison method,  we will show that  the expectation value of this term will obtain an extra factor $N^{-1/2}$, i.e., 
\be\label{BN3}
 |\E\,  h(t_X)(Y G)_{ij}|\leq CN^{-1/2+\sigma}\sqrt{\frac{|w|^{1/2}}{N\eta} }
\ee
This extra factor $N^{-1/2}$ plays a key role in our new proof.  A similar method was   used in the \cite{BEKYY2012}.  

Now we explain the  basic idea of proving \eqref{BN3}-type bounds, i.e.,where the extra $N^{-1/2}$ factor comes from.  For simplicity we assume $X_{ij}\in\R$. Let $Y^{(i,i)}_z$ be the matrix obtained by removing  $i-$th row and column of $Y_z$, and define 
$$G^{(i,i)}:=((Y^{(i,i)}_z)^*Y^{(i,i)}_z-w)^{-1}, \quad
 \mG^{(i,i)}:=(Y^{(i,i)}_z(Y^{(i,i)}_z)^*-w)^{-1}
 $$
We  write $h(t_X)(Y_zG)_{ij}$ as the polynomials of the $i-$th row/column of $X$: $X_{ik}$, $X_{ki}$ ($1\leq k\leq N$), $G^{(i,i)}$ and $ \mG^{(i,i)}$, i.e., 
$$
h(t_X)(Y_zG)_{ij}=P(\{X_{ik}\}_{k=1}^N,\; \{X_{ki}\}_{k=1}^N,\; G^{(i,i)}, \, \mG^{(i,i)})+{\rm negligible\; error }
$$
 where $P$ is a polynomial.   By definition,  $X_{ik}$, $X_{ki}$ are independent of $G^{(i,i)}$ and $\mG^{(i,i)}$.  In this polynomial, we will show that the degrees of every monomials w.r.t. $X_{ik}$ and $X_{ki}$'s are always odd numbers.  Therefore, in taking the expectation value, with assumption $\E X_{ij}=0$ and $|\E X^k_{ij}|\leq O(N^{-k/2})$, one will see an extra combination factor $N^{-1/2}$. The following simple example will show why the odd powers give an extra factor $N^{-1/2}$. Suppose we estimate $\E\sum_{kst} X_{ik}G^{(i,i)}_{kl}X_{is}G^{(i,i)}_{st}X_{it}$.  Since $X_{ik}$, $X_{ki}$ are independent of $G^{(i,i)}$ and $\mG^{(i,i)}$, $\E X_{ij}=0$ and $|\E X^k_{ij}|\leq O(N^{-k/2})$, the nonzero contributions only come  from the terms where  $k=s=t$, therefore
$$
| \E \sum_{kst} X_{ik}G^{(i,i)}_{kl}X_{is}G^{(i,i)}_{st}X_{ti}| = |\E \sum_{k } X_{ik}G^{(i,i)}_{kl}X_{ik}G^{(i,i)}_{st}X_{ki}|\le C N^{-1/2}\E(\max_{ab}|G^{(i,i)}_{ab}|)^2 
$$
On the other hand, without $\E $, this term can only be bounded without this $N^{-1/2}$ factor (with large deviation theory).
 $$
|   \sum_{kst} X_{ik}G^{(i,i)}_{kl}X_{is}G^{(i,i)}_{st}X_{ti}|=|   \sum_{k } X_{ik}G^{(i,i)}_{kl}||\sum_{st}X_{is}G^{(i,i)}_{st}X_{ti}|\leq (\log N)^C\, (\max_{ab}|G^{(i,i)}_{ab}|)^2 
$$
Note: one will not see this $N^{-1/2}$ factor if the degree is even number, e.g., $\E\sum X_{is}G^{(i,i)}_{st}X_{it}$.  Based on this new idea, the main  task  of proving  lemma \ref{main lemma} and \eqref{BN3}-type bounds is writing the functionals of $Y_z$'s, $G$'s and $\mG$'s as  the polynomials of $X_{ik}$, $X_{ki}$ ($1\leq k\leq N$) ,  $G^{(i,i)}$ and  $\mG^{(i,i)}$ for some $1\leq i\leq N$,  
(up to negligible error)  and counting the degree of each monomial.

 \section{Proof of Lemma \ref{main lemma}}
 
In this section, we apply the Green's function comparison method to prove the Lemma \ref{main lemma}. We will see the key input of  proving  Lemma \ref{main lemma} is the  lemma \ref{mainnewL}.  This new lemma is similar to (3.62)-(3.63) of \cite{BouYauYin2012Edge}, but without the third moments vanishing assumption.  More precisely,  the (3.62)-(3.63) of \cite{BouYauYin2012Edge} is similar to the \eqref{36ol} of this work, and lemma \ref{mainnewL} is the key step of proving \eqref{36ol}.   The proof of lemma \ref{mainnewL} will start from section \ref{sec: 4}. In \cite{BouYauYin2012Edge}, the (3.62)-(3.63)  can be easily proved by bounding 
 the expectation value of these terms with their stochastically dominations.  In this paper, as introduced in subsection \ref{subsection: bs}, we will introduce a new comparison method to show that, for the contribution comes from $X_{ij}$'s third moment,     their expectation values    have an extra factor $N^{-1/2}$, i.e., lemma \ref{mainnewL}. 
 
  First of all,  we state the following lemma. It will be used to estimate the expectation value of some random variables which are stochastically dominated, but not $L_\infty$ bounded.

\begin{lemma}\label{dzzy}
Let $v=v^{(N)}$ be a family of centered random variables with variance $1/N$, satisfying the sub exponential decay \eqref{subexp}.  
Let $\wt A=\wt A^{(N)}$ and $A=A^{(N)}$ be   families of   random variables. Suppose   $A\prec 1$, and $A=\sum_{n=0}^{C} A_n\,v^n$, where $|A_n|\le N^C$ for some fixed constant $C>0$. We also assume that $\wt A$ is independent of $v$ and $|\wt A|\leq N^C$ for some $C>0$. Then for any fixed $ p\in \N$ and fixed (small) $\delta>0$, 
\be\nonumber
|\E\,\wt A A \,v^p| \le (\E |\wt A|  )N^{-p/2+\delta}+N^{-1/\delta}
\ee
 for large enough $N$. 
\end{lemma}
Note: Here   $A$ or $A_i$'s may depend on $v$.

\bigskip

  {\it Proof of Lemma \ref{dzzy}: } By definition \ref{stodom}, the assumption $A\prec 1$, and the fact that $v$ has sub exponential decay \eqref{subexp}, for any fixed $\delta>0$ and $D>0$  there is a probability subset $\Omega$ such that $\P(\Omega)\ge 1-N^{-D}$ and 
$$
|1_\Omega Av^p|\leq  N^{-p/2+\delta}
$$
Then 
\be\nonumber
|\E\wt A A v^p|\leq  (\E |\wt  A|) N^{-p/2+\delta}+|\E_{\Omega^c} \wt AA v^p|\leq (\E |\wt  A|)  N^{-p/2+\delta}+O( N^{-D/2+2C})
\ee
for the second inequality, we used Cauchy Schwarz inequality. Choosing large enough $D$, we complete the proof of lemma \ref{dzzy}. 

\qed

 Because of this lemma, for any centered random variables $v$ with variance $1/N$, satisfying the sub exponential decay \eqref{subexp}, we define
\be\label{defMv}
\cal M_C(v):= \left\{A: A=\sum_{n=0}^{C} A_n v^n, \;|A_n|\leq N^C  
 \right\}
\ee
Now we return to prove Lemma \ref{main lemma}.

{\it Proof of Lemma \ref{main lemma}: }  For simplicity, we assume that the matrix entries are real numbers. Let $X$ and $X'$ be two ensembles which satisfy the assumption of Theorem \ref{z1}. 
To prove Lemma \ref{main lemma}, as we explained in the beginning of subsection 2.1 (near \eqref{Bn2}), one only needs to show that for any  fixed small enough $\e>0$, $ s\in (0,1/2)$, and $p\in 2\N$, if $||z_0|-1|\leq \e$  then
 \be \label{scon}
 \E |A_{X}^{(f)}|^p  \leq C\,\E |A_{X'}^{(f)}|^p+N^{C\e p},  
\ee
for large enough $N$.  For integer $k$, $0\leq k\leq N^2$, define the following matrix $X_k$ interpolating between $X'$ and $X $:
$$
X_k(i,j)=
\left\{
\begin{array}{ccc}
X (i,j)&{\rm\ if\ }&k\ge N(i-1)+j\\
X' (i,j)&{\rm\ if\ }&k<N(i-1)+j\\
\end{array}
\right..
$$
Note that $X'=X_0$ and $X =X_{N^2}$.  As one can see that the difference between $X_k$ and $X_{k-1}$ is just one matrix entry. We denote the index of  this entry as $(a,b):=(a_k, b_k)$ ($a_k$, $b_k\in \Z$,   $1\leq a_k, b_k\leq N$), here  $k=(a_k-1)N+b_k$

 Furthermore, we define $t_{X_{k-1}}$,   $t_{X_{k}}$,  $A^{(f)}_{X_{k-1} }$,   $A^{(f)}_{X_{k} }$ with $X_{k-1}$ and $X_k$, as in Def. \ref{DefAX}. 
 We are going to show that if this special matrix entry is in the diagonal line, i.e., $a=b$ then 
 \be\label{35ol}
 \left | \E\left(A^{(f)}_{X_{k } } \right)^p-\E\left(A^{(f)}_{X_{k-1 } } \right)^p \right | \le N^{-3/2} \left( N^{\e}+2\E\left( A^{(f)}_{X_{k-1} }\right) ^p\right)
 \ee
 otherwise, i.e., $a\neq b$,  
 \be\label{36ol}
 \left | \E\left(A^{(f)}_{X_{k } } \right)^p-\E\left(A^{(f)}_{X_{k-1 } } \right)^p \right | \le N^{-2} \left( N^{\e}+2\E\left( A^{(f)}_{X_{k-1} }\right) ^p\right)
 \ee
for sufficiently large $N$ (independent of $k$). Clearly,  \eqref{35ol} and \eqref{36ol} imply \eqref{scon}. 

 We are going to compare the these functionals corresponding to $X_k$ and $X_{k-1}$ with a third one,
corresponding to the matrix $\wt Q$ hereafter with deterministic $(a,b)$ entry. We define
the following $N\times N$ matrices (hereafter, $Y_\ell=X_\ell-z I$,  $\ell= k$ or $k-1$):
\begin{align}
v&=v_{ab}{\bf e}_{ab}=X'(a,b){\bf e}_{ab},\label{def:vab}\\
u&=u_{ab}{\bf e}_{ab}=X(a,b){\bf e}_{ab},\label{def:uab}\\
\wt Q&=X_{k-1}-v=X_{k}-u,\label{def:wtQ}\\
Q&=Y_{k-1}-v=Y_{k}-u,\label{def:Q}\\
R&=(Q^* Q-wI)^{-1}\label{def:R}\\
\cal R&=(Q Q^*-wI)^{-1}\label{def:calR}\\
S&=(Y^*_{k-1} Y_{k-1}-w I)^{-1}\label{def:S}\\
T&=(Y^*_{k} Y_{k}-w I)^{-1}\label{def:T}
 \end{align}

Furthermore, we define  $t_{\wt Q}$,  $A^{(f)}_{\wt Q }$  with $\wt Q$, as in Def. \ref{DefAX}. To prove \eqref{35ol} and  \eqref{36ol}, we will estimate
$  A^{(f)}_{X_{k-1 } } - A^{(f)}_{X_{\wt Q} } $ and $  A^{(f)}_{X_{k  } } - A^{(f)}_{X_{\wt Q} } $ .  

\medskip

First we introduce the notations 
$$
m_S=\frac{1}{N}\tr S, \quad m_R=\frac{1}{N}\tr R, \quad m_T=\frac{1}{N}\tr T
$$
We note: with   Cauchy's interlace theorem,   for some $C>0$,
\be\label{iwhmm}
 |m_S-m_R|
 \leq C(N\eta)^{-1}, \quad \eta=\im w
\ee
 holds  for any $w$ and $z$.  It implies 
\be\label{yongd}
 |A^{(f)}_{X_{k-1}}-A^{(f)}_{\wt Q }|
 \leq C
\ee
To estimate $A^{(f)}_{X_{k-1}}-A^{(f)}_{\wt Q }$,  from \eqref{defAXf},  we have 
\be\label{sdsn}
A^{(f)}_{X_{k-1}}-A^{(f)}_{\wt Q }= N
\int  \Delta f(\xi)  \int_{I}  \chi(\eta)\phi '(E)
\Big(h(t_{X_{k-1}})\re m_S-h(t_{\wt Q})\re m_R\Big) \rd E\rd \eta\rdA(\xi), \quad  
 \ee
where $z=z_0+N^{-s}\xi$, $w=E+i\eta$,  $\phi=\phi_{\e,z}$.  Recall $t_{X_{k-1}}$ and $t_{\wt Q}$ are defined with $m_S$ and $m_R$ respectively.  Applying Taylor's expansion on the term $ h(t_{X_{k-1}})\re m_S-h(t_{\wt Q})\re m_R $ in \eqref{sdsn} and letting $h^{(k)}$ be the $k$th derivative of $h$, we have   
\begin{align}\label{ygldj}
h(t_{X_{k-1}})\re m_S-h(t_{\wt Q})\re m_R=& \sum_{ n=1}^3B_n(\wt Q)\left(\re m_S-\re m_R\right)^ n+B_4(X_{k-1}, \wt Q)\left(\re m_S-\re m_R\right)^4
 \end{align}
where   $B_n(\wt Q)$ ($1\leq n\le 3$)  and $B_4(X_{k-1}, \wt Q) $ are  defined as 
\begin{align}\label{defBnn}
B_n(\wt Q):= \frac1{n!}(N^{1-\e}\eta )^{ (n-1)}\left(nh^{( n-1)}(t_{\wt Q})+h^{( n)}(t_{\wt Q})t_{\wt Q}\right)
\\\nonumber
B_4(X_{k-1}, \wt Q):=\frac{1}{24}(N^{1-\e}\eta)^3 \left(4h^{(3)}(\zeta) +h^{(4)}(\zeta)  \zeta\right)
\end{align}
where  
 $\zeta$ is between $t_{X_{k-1}}$ and $t_{\wt Q}$, and only depends on $t_{X_{k-1}}$, $t_{\wt Q}$ and $h$. As one can see that $B_1$, $B_2$ and $B_3$ are independent of $v_{ab}$. For  the definition of $B $'s, we note that if $n\ge 1$, then 
  $$  
h^{(n)}(x)\neq 0\implies x\sim1
$$
Therefore, with $|h|\leq 1$, we obtain the following uniform bounds for $B$'s:
\be\label{Bbound}
|B_n|\leq  (N^{1-\e}\eta)^{(n-1)}, \quad 1\leq n\leq 4
\ee
To estimate the $m_S-m_R$ in \eqref{ygldj}, we study the difference between $m_S$ and $m_R$ in the parameter set:
\be\label{paset}
\left\{(k, z, w)\in \mathbb Z\times \C^2\;:\; 0\leq k\leq N^2, \; ||z|-1|\leq2\e, \;w\in I_\e\right\}
\ee
Recall in  (3.59) of [2] and the discussion below (3.61) of [2], it was  proved that 
with the notations:
\begin{align}\label{PPP}
&P_1 ( \wt Q) :=\frac1N \re \left(-  2( QR^2)_{ab} \right)
\\\nonumber
 &P_2 ( \wt Q) :=\frac1N \re \left(w\cal R_{aa} (R^2)_{bb}+  2  ( QR^2)_{ab} ( R  Q^*)_{ba} + (QR^2Q^*)_{aa} R_{bb}\right)
\\\nonumber
 &P_3 ( \wt Q) :=
\\\nonumber
&\frac1N \re \left(- 2 ( R  Q^*)^2_{ba}   (QR^2)_{ab}  -2( R  Q^*)_{ba} (QR^2 Q^*)_{aa} R_{bb}
- 2( R  Q^*)_{ba}w\cal R_{aa} (R^2)_{bb} 
-2 w\cal R_{aa}R_{bb}  (QR^2)_{ab} )\right)
 \end{align}
  the difference between $\re m_S$ and $\re m_R$, i.e., $(\frac1N\re\tr S-\frac1N\re\tr R)$ can be written  as (recall $v_{ab}=X'(a,b)$) 
 \be\label{lryy}
\re m_S-\re m_R=\sum_{n=1}^3 P_n(\wt Q) \cdot (v_{ab}) ^3+  P_4(X_{k-1}, \wt Q) \cdot(v_{ab})^4,
\ee
where $P_4(X_{k-1}, \wt Q)$ depends on $X_{k-1}$ and $\wt Q$, and  the $P$'s can be bounded as 
\be\label{PPboun}  P_1(\wt Q),\;  P_2(\wt Q), \;P_3(\wt Q), \; P_4(X_{k-1}, \wt Q)\prec  (N\eta)^{-1},  \ee 
uniformly  for 
  $(k, z, w)$ in \eqref{paset}.  In [2], the uniformness was not emphasized, but it can be easily checked.   From \eqref{def:wtQ}-\eqref{def:calR} and the definition of $P_{1,2,3}(\wt Q)$, we can see that $P_{1,2,3}(\wt Q)$ only  depend  on $\wt Q$ and they are independent of $v_{ab}$.

Now we collect some simple bounds on $P_i$'s. For $L_\infty$ norm, by definition, 
 it is easy to prove that the following inequalities always hold: 
 \be\nonumber
\|S\|,\quad  \|\cal R\|, \quad  \|R\|, \quad  \|R^2\|, \quad \| QR\|, \quad\|  QR^2\|, \quad\| QR^2 Q^*\|\le  N^C
 \ee 
 for 
 any $(k, z, w)$ in \eqref{paset} and some fixed constant $C>0$. Then with the definition in \eqref{PPP},  we also have that for any $(k, z, w)$ in \eqref{paset} and some constant $C>0$
   \be\label{PPPNC}
 P_1, P_2, P_3 =O(N^C).
\ee
  Expanding $S$ around $R$ with the fact: $S=  (R^{-1}+(Y_k^*Y_k- Q^* Q) )^{-1}$,  we obtain that  for any fixed $m\in \mathbb N$
\be\label{317s}
S-R=\sum_{n=1}^m\left(-R(Y_k^*Y_k- Q^* Q)\right)^nR+\left(-R(Y_k^*Y_k- Q^* Q)\right)^{m+1}S
\ee
Let $m=5$ in \eqref{317s}.  Now we take  $\frac{1}{N}\re \tr  $ on the both sides of \eqref{317s} and compare it with \eqref{lryy}. Since
$Y_k^*Y_k- Q^* Q=v _{ab}({\bf e}_{ba}Q)+v_{ab}(Q^*{\bf e}_{ab})+v^2_{ab}{\bf e}_{aa}$,  we can  see that for $1\leq l\leq 3$, the $P_l(Q)$ is the coefficient of the $(v_{ab})^l$ term in the r.h.s. of $\frac{1}{N}\re \tr  \eqref{317s}$ and 
\be\label{cmdq}
P_4(X_{k-1}, \wt Q)\in \cal M_C(v_{ab}) {\cor }
\ee 

Similarly, using this expansion ($m=5$), and the fact:
$$
\partial_w R= R^2=O(N^C),\quad \partial _z R= R(Q+Q^*)R=O(N^C),
$$
and $\partial_w S$, $\partial_z S=O(N^C)$,  we can improve \eqref{PPboun} to the following one: 
\be 
\max _{(k, z, w)\in \eqref{paset}}\label{PPboun2}  N\eta\left(|P_1(\wt Q)|+|  P_2(\wt Q)|+|P_3(\wt Q)|+| P_4(X_{k-1}, \wt Q)  |\right)\prec 1
 \ee 
 We note: this statement shows that \eqref{PPboun} can hold for different $ (k, z, w ) \in \eqref{paset}$ with the same probability subset.

Inserting \eqref{ygldj} and \eqref{lryy} into   \eqref{sdsn}, we   write $A^{(f)}_{X_{k-1} }-A^{(f)}_{Q }$ as a polynomial of $v_{ab}$ as follows.  
 \begin{align}\label{eqn:recY}
A^{(f)}_{X_{k-1} }- A^{(f)}_{\wt Q } =\htP_1(\wt Q) \cdot v_{ab}+ \htP_2(\wt Q)  \cdot (v_{ab})^2+\htP_3(\wt Q)  \cdot (v_{ab})^3+\htP_4(X_{k-1}, \wt Q)    \cdot (v_{ab})^4.
\end{align}
where \begin{align}\label{YYY}
\htP_1 (\wt Q): =&N
\int  \Delta f(\xi)  \int_{I} \big(B_1P_1 \big) \chi(\eta)\phi '(E)
  \rd E\rd \eta\rdA(\xi)\\\nonumber
\htP_2 (\wt Q):  =&N
\int  \Delta f(\xi)  \int_{I}  \big(B_1P_2+ B_2P_1^2  \big) \chi(\eta)\phi '(E)
  \rd E\rd \eta\rdA(\xi)\\\nonumber
\htP_3  (\wt Q): =&N
\int  \Delta f(\xi)  \int_{I}  \Big( B_1P_3+2B_2P_1P_2+B_3P_1^3\Big)  \chi(\eta)\phi '(E)  \rd E\rd \eta\rdA(\xi) \\\nonumber
\htP_4  (X_{k-1}, \wt Q): =&N
\int  \Delta f(\xi)  \int_{I}  \Big( \sum_n B_n \sum _{\sum_j i_j\ge 4} (v_{ab})^{(\sum_j i_j)-4} \prod_{j=1}^n P_{i_j}  \Big)  \chi(\eta)\phi '(E)  \rd E\rd \eta\rdA(\xi) 
 \end{align}
where $B_n=B_n({\wt Q})$, $P_n=P_n(\wt Q)$ ($1\le n\le 3$), $B_4=B_4(X_{k-1}, \wt Q)$ and $P_4=P_4(X_{k-1}, \wt Q)$.  
We note: 
$ \htP_1(\wt Q)$, $\htP_2(\wt Q)$ and $\htP_3(\wt Q)$ are independent of $v_{ab}$.  

Replacing $X_{k-1}$ with $X_{k}$, with the same method, we obtain (Here $v_{ab}$ is replaced with $u_{ab}$)
\begin{align}\label{eqn:recY2}
A^{(f)}_{X_{k} }- A^{(f)}_{\wt Q } =\htP_1(\wt Q) u_{ab}+ \htP_2(\wt Q) u^2_{ab}+\htP_3(\wt Q) u^3_{ab}+\htP_4(X_{k}, \wt Q)    u_{ab}^4 
\end{align}

From \eqref{Bbound} and \eqref{PPboun2}, it is easy to check that $\htP_1$, $\htP_2$ and $\htP_3\prec 1$ uniformly hold for $1\leq k\leq N^2$. For $L_\infty$ bound, with \eqref{PPPNC},  they are bounded by $N^C$ for some $C$.  Similarly, we can obtain that 
$\htP_4\prec 1$. With 
\eqref{cmdq}, we have $\htP_4  (X_{k-1}, \wt Q)\in \cal M_C(v_{ab})$ and $\htP_4  (X_{k}, \wt Q)\in \cal M_C(u_{ab})$. So far, we proved 
 \be\label{thpbou}
 \htP_{1,2,3,4}\prec 1, \quad  \htP_{1,2,3}(\wt Q)\leq N^C,    \quad \htP_4  (X_{k-1}, \wt Q)\in \cal M _C(v_{ab}), \quad \htP_4  (X_{k}, \wt Q)\in \cal M_C(u_{ab})
 \ee
 uniformly hold for $1\leq k\leq N^2$. 
 
 Now we return to  prove \eqref{35ol} and \eqref{36ol}. First we write  $$
(A^{(f)}_{X_{k-1} })^p-(A^{(f)}_{X_{k} })^p=\sum_{j=0}^{p-1} \binom{p}{j}
(A^{(f)}_{\tilde Q })^j \left((A^{(f)}_{X_{k-1} }- A^{(f)}_{\wt Q }) ^{p-j}-(A^{(f)}_{X_{k} }- A^{(f)}_{\wt Q })^{p-j}\right).
$$

 We insert the  \eqref{eqn:recY} and \eqref{eqn:recY2} into the r.h.s. and write it in the following form
 \be\label{AmBm}
 (A^{(f)}_{X_{k-1} })^p-(A^{(f)}_{X_{k} })^p=\sum_{m=1}^{4p} \left(\cal A_{m} v^m_{ab}- \cal B_{m} u^m_{ab}\right)\ee
where $\cal A_{m}$ only contains $A^{(f)}_{\tilde Q }$, $\htP_{1,2,3}(\wt Q)$, $\htP_4  (X_{k-1}, \wt Q)$, and $\cal B_{m}$ only contains $A^{(f)}_{\tilde Q }$, $\htP_{1,2,3}(\wt Q)$, $\htP_4  (X_{k}, \wt Q)$  For example,   
\be\label{defcalA}
\cal A_{3}=\cal B_{3}=C_{p,3}(A^{(f)}_{\wt Q })^{p-3}\htP_1^3(\wt Q)
+C_{p,2}(A^{(f)}_{\wt Q })^{p-2}\htP_1(\wt Q)\htP_2(\wt Q)+C_{p,1}(A^{(f)}_{\wt Q })^{p-1}\htP_3(\wt Q)
\ee
where $C_{p,n}$ ($1\leq n\le 3$) are constants only depends on $p$.  Since   the first two moments of $v_{ab}$ and $u_{ab}$ coincide, $u_{ab},v_{ab}$ are independent of  $ \wt Q$, and $\cal A_{1}=\cal B_1$, $\cal A_{2}=\cal B_2$ only contain $A^{(f)}_{\tilde Q }$, $\htP_{1,2,3}(\wt Q)$, we have 
$$
\E  (A^{(f)}_{X_{k-1} })^p-\E (A^{(f)}_{X_{k} })^p=\sum_{m=3}^{4p}  \left(\cal A_{m} v^m_{ab}- \cal B_{m} u^m_{ab}\right) $$ 
Recall the definition of $\cal A_m$ and $\cal B_m$ from \eqref{AmBm}, for the terms $m\ge 4$, using \eqref{thpbou} and Lemma \ref{dzzy}, we get
$$
|\E\sum_{m=4}^{4p}  \left(\cal A_{m} v^m_{ab}- \cal B_{m} u^m_{ab}\right) |\leq \sum_{j=0}^{p-1}
\E\left|(A^{(f)}_{\wt Q })^j\right| \OO_\prec(N^{-2})+N^{-2} 
\le
N^{-2}\left( \OO_\prec(1)+\E |A^{(f)}_{\wt Q }|^p \right).
$$
Therefore, with $\cal A_3=\cal B_3$, 
\begin{equation}\label{eqn:rec2}
\left | \E\left(A^{(f)}_{X_{k-1} })\right)^p-\E\left(A^{(f)}_{X_{k } })\right)^p \right | \le
N^{-2}\left( \OO_\prec(1)+\E |A^{(f)}_{\wt Q }|^p \right)+|\E \cal A_3| \left(|\E v^3_{ab}|+|\E u^3_{ab}|\right),
\end{equation}
Similarly,  using \eqref{thpbou}, \eqref{defcalA}, $A_{\wt Q}^{(f)} =O( N^C)$ and Lemma \ref{dzzy},  we
 have
 \be\label{ECAw}
|\E \cal A_3| \left(|\E v^3_{ab}|+|\E u^3_{ab}|\right)\leq N^{-3/2}\left( \OO_\prec(1)+\E |A^{(f)}_{\wt Q }|^p  \right)
 \ee
 As in [2]-(3.64), using H\"older's inequality and the bound \eqref{yongd}, we have
\begin{align}\label{bhjq}
\E |A^{(f)}_{\wt Q }|^p&\leq \E\left((A^{(f)}_{X_{k-1} })^p\right)+
\sum_{j=1}^p\binom{j}{p}\E\left(\left|A^{(f)}_{X_{k-1} }\right|^{p-j}|A^{(f)}_{X_{k-1} }- A^{(f)}_{\wt Q }|^j\right)\\\nonumber
&\leq
\E\left((A^{(f)}_{X_{k-1} })^p\right)+
\sum_{j=1}^p\binom{j}{p}\E\left(\left|A^{(f)}_{X_{k-1} }\right|^{p}\right)^{\frac{p-j}{p}}
\E\left(\left|A^{(f)}_{X_{k-1} }- A^{(f)}_{\wt Q }\right|^{p}\right)^{\frac{j}{p}}, \\\nonumber
&\leq \left(\OO_\prec(1)+2\E |A^{(f)}_{X_{k-1}}|^p\right),
\end{align}
Then combining \eqref{eqn:rec2}-\eqref{bhjq},  we obtain \eqref{35ol}. (Note: $p\in 2\mathbb Z$.)
 
To prove \eqref{36ol},  we claim the following lemma, which provides the stronger bound  on  the expectation value of  the r.h.s. of \eqref{defcalA}.

 \begin{lemma}\label{mainnewL} Assume $1\leq a\neq b\leq N$. Let $X$ be defined as in  Theorem \ref{z1},  except that $X_{ab} =0$. For any fixed small
enough $\e>0$, if $||z_0|-1|\leq \e$ and $s\in(0,1/2)$, define $A^{(f)}_{X }$, $P_i(X)$, $B_i(X)$, $\htP_i(X)$, $i=1,2,3$ as in \eqref{defAXf}, \eqref{PPP}, \eqref{defBnn} and \eqref{YYY}. (More precisely $\wt Q$, $Q$, $R$ in  \eqref{PPP} and  \eqref{defBnn} will be replaced with $X$, $Y=X-zI$ and $(Y^*Y-wI)$ respectively.)  Then
\be\label{110X}
\left|\E(A^{(f)}_{X })^{p-3}\htP_1^3(X)\right|+\left| \E (A^{(f)}_{X })^{p-2}\htP_1(X)\htP_2(X)\right|+
\left| \E  (A^{(f)}_{X })^{p-1}\htP_3(X)\right|\prec N^{-1/2}
\left( \OO_\prec(1)+\E |A^{(f)}_{X}|^p \right)
\ee
uniformly for $(a,b)$.
\end{lemma}

We return to prove \eqref{eqn:rec2} and prove lemma \eqref{mainnewL} in the next section. Inserting this lemma and \eqref{defcalA} into \eqref{eqn:rec2}, as in \eqref{ECAw},  we obtain   that if  $a\neq b$, then
\be\label{ECAw2}
 |\E \cal A|\leq N^{-2}\left( \OO_\prec(1)+\E |A^{(f)}_{\wt Q}|^p\right)
 \ee
Together with \eqref{eqn:rec2} and \eqref{bhjq}, we   obtain \eqref{36ol}.  Clearly,  \eqref{35ol} and \eqref{36ol} imply \eqref{scon}, and we complete the proof of lemma \ref{main lemma} and lemma \ref{lyy}. 

\qed

 \section{Proof of Lemma \ref{mainnewL}}\label{sec: 4}

 Lemma \ref{mainnewL} bounds the expectation values of some  polynomials of $A^{(f)}_{X}$ and $\htP_{1,2,3} (X)$. Roughly speaking   Lemma \ref{mainnewL} shows that the expectation value of these polynomials are much less than their stochastic domination by a factor $N^{-1/2}$. (Note: $a$ and $b$ appear in the definitions of $P_{1,2,3}$ and $B_{1,2,3}$. The $\htP_{1,2,3}$ are defined with $P_{1,2,3}$ and $B_{1,2,3}$.) As introduced in the second part of subsection \ref{subsection: bs} (below \eqref{BN3}), the main strategy of showing this extra factor  is 
 \begin{itemize}
\item writing them as the polynomials  (up to negligible error) of   $X_{ak}$'s, $X_{ka}$'s ($1\leq k\leq N$),   
 $G^{(a,a)}$ and  $\mG^{(a,a)}$, which are defined as 
 $$G^{  (a,a)}:=((Y^{  (a,a)}_z)^*Y^{  (a,a)}_z-w)^{-1}, \quad
 \mG^{  (a,a)}:=(Y^{  (a,a)}_z(Y^{  (a,a)}_z)^*-w)^{-1}
 $$
 and $Y^{(a,a)}:=Y^{(a,a)}_z$ is  the matrix    obtained by removing  the $i-$th row and column of $Y_z$.
 \item showing the degrees of the monomials of $X_{ak}$'s and $X_{ka}$'s in above polynomials  are always odd (except for $X_{aa}$).  
\end{itemize}

First of all, in lemma \ref{exp22} and \ref{exp223}, we introduce  some polynomials having the properties we need for Lemma \ref{mainnewL}, i.e., their expectation values have an extra factor $N^{-1/2}$ comparing with their stochastic  domination. In the next subsection,  we introduce some $\cal F$ sets, whose elements are the   "basic" polynomials in our proof, i.e., the bricks of the polynomials in lemma \ref{exp22} and \ref{exp223}.

 \subsection{Basic polynomials and their properties. }

We first introduce some notations.

\begin{definition} \label{definition of minor} {$X^{(\bb T, \bb U)}$, $Y^{(\bb T, \bb U)}$, $G^{(\bb T, \bb U)}$
 and $\mG^{(\bb T, \bb U)}$}

Let $\bb T, \bb U$ be some subsets of $\{1,2,\cdots, N\}$. Then we define $Y^{(\bb T, \bb U)}$ as the $ (N-|\bb U|)\times (  N-|\bb T|)  $
matrix obtained by removing all columns of $Y$   indexed by $i \in \bb T$ and all rows  of $Y$  indexed by $i \in \bb U$. Notice that we keep the labels of
indices of $Y$ when defining $Y^{(\bb T, \bb U)}$. With the same method, we define $X^{(\bb T, \bb U)}$ with $X$.

Let $\by_i$ be the $i $-th column of $Y$ and  $\by^{(\bb S)}_i$ be the vector obtained by removing $\by_i (j) $ for
 all  $ j \in  \bb S$. Similarly we define $\mathrm y_i$ be the $i $-th row of $Y$.
Define
\begin{align*}
G^{(\bb T, \bb U)}=\Big [  (Y^{(\bb T, \bb U)})^* Y^{(\bb T, \bb U)}- w\Big]^{-1},\ \  & m_G^{(\bb T, \bb U)} =\frac{1}{N}\tr G^{(\bb T, \bb U)},
\\
\mG^{(\bb T, \bb U)}= \Big [ Y^{(\bb T, \bb U)}(Y^{(\bb T, \bb U)})^*- w \Big ]^{-1},\ \ & m_\mG^{(\bb T, \bb U)} =\frac1N\tr \mG^{(\bb T, \bb U)}.
\end{align*}
By definition,  $m^{(\emptyset, \emptyset)} = m$.
 Since the eigenvalues of $Y^* Y $ and $Y Y^*$ are the same except the zero eigenvalue, it is easy to check that
\be\label{35bd}
m_G^{(\bb T, \bb U)}(w) =m_\mG^{(\bb T, \bb U)} +\frac{|\bb U|-|\bb T|}{N  w}
\ee
For $|\bb U|=| \bb T|$, we define
\be\label{d trGmG}
m ^{(\bb T, \bb U)}:= m_G^{(\bb T, \bb U)} = m_\mG^{(\bb T, \bb U)}
\ee
\end{definition}

There is a crude bound for $(m_G^{(\bb T, \bb U)}-m)$ proved in (6.6) of \cite{BouYauYin2012Bulk}:
\be\label{boundmGm}
\left|m_G^{(\bb T, \bb U)}-m\right|+\left|m_\mG^{(\bb T, \bb U)}-m\right|\leq C\frac{|\bb T|+|\bb U|}{N\eta}
\ee
 \begin{definition}\label{def: Set}{\bf Notations for general sets.}
 
  As usual, if  $x\in \R$ or $\C$,  and $\cal S$ is a set of random variables then $x \cal S$ denotes the following set  as
 $$
x\cal S:= \left\{ x s: s\in \cal S\right\}
 $$
 For two sets $\cal S_1$ and $\cal S_2$ of random variables, we define  the following set  as
 $$
 \cal S_1\cdot \cal S_2:=\left\{s_1\cdot s_2\;\big|\; s_1\in \cal S_1, \; s_2\in \cal S_2\right\}
 $$
 For simplicity, we call $s\in_n \cal S$ if and only if   $s$ can be written as the sum of $O(1)$ elements in $\cal S$, i.e., 
\be\nonumber
s\in_n \cal S\iff s\in \left\{\sum_{i=1}^n s_i\Big|\; s_i\in \cal S, \;n\in \N, \;  n=O(1)\right\}
\ee
  \end{definition}

 \begin{definition}\label{def: FFF} {\bf Definition of $\cal F_0$, $\cal F_1$, $\cal F_{1/2}$ and  $\cal F $.} 
 
 For fixed indeces $a,b$ and ensemble $X$ in lemma \ref{mainnewL}, we define $\cal F_0$ as the set of random variables (depending on $X$) which are  stochastically dominated by 1 and independent of any $X_{ak}$ and $X_{ka}$ $(1\leq k\leq N)$,  i.e., 
 \begin{align}\nonumber
  \cal F_0 =\left\{  V\;:\;   V\prec 1, \;   V \; {\rm is \; independent \;of \;
 the} \; a{\rm -th\; row \; and \; column \;of\; X  } \right\}
 \end{align}
Note: $\cal F_0$ depends on $a$, not $b$. One example element in $\cal F_0$ is $\tr X-X_{aa}$.

 \medskip
 
 For simplicity, we define
 $$
 \sum _{i}^{(a)}:= \sum_{i\neq a}, \quad\quad  \sum _{ij}^{(a)}:= \sum_{ij\neq a}
 $$
 
Next we define $\cal F_1$ as the union of the set $( N^{ 1/2}X_{aa} \cal F_0)$ and  the sets of some quadratic forms as follows
 \begin{align}\nonumber
\cal F_1:= &
\left( N^{ 1/2}X_{aa} \cal F_0  \right)\bigcup\left\{  
 \sum_{kl}^{(a)}X _{ka}V_{kl}X_{la}  \;{\rm or }\;\sum_{kl}^{(a)}X _{ak}V_{kl}X_{al} 
 \;\Bigg| \; \max_{k  l}|V_{kl}|\prec 1,\; V_{kl} \; \in \cal F_0  \right\}
\\\nonumber  
  \bigcup & \left\{  \sum_{k\neq l}^ {(a)}X_{ak}V_{kl}X_{la} +N^{1/2}\sum_{k }^{(a)}X_{ak}V_{kk}X_{ka} 
 \;\Bigg|\;\max_{k  l}|V_{kl}|\prec 1,   \;V_{kl} \in \cal F_0  \right\}
 \end{align}
(Note it is $X _{ka}V_{kl}X_{la}$ or $X _{ak}V_{kl}X_{al}$ in the first line and $X_{ak}V_{kl}X_{la}$ in the second line, and the diagonal terms in the second case is allowed to be larger than the others by a factor $N^{1/2}$.) 

\bigskip

Furthermore, we define  $\cal F$ as the set of following random variables  
\be\nonumber
\cal F :=\left\{V\;\Big|\;V\in_n\cal F_0 \bigcup \left( \bigcup_{n=O(1)}\,
  \left(\cal F_1 \right)^n \right)\; \right\} \;  
\ee
 where $\left(\cal F_1\right)^n$ represents the set of  the  products of $n$ elements in $\cal F_1$.  For simplicity, sometimes we write $\cal F$
 $$
 \cal F=\cal F_\emptyset
 $$
  i.e., with the subscription empty set $\emptyset$.

  Similarly, we define 
 \be\nonumber
 \cal F_{1/2}  =\left\{  \sum_{k }^{(a)}X_{ak}V_{k } \;{\rm or} \;  \sum_{k }^{(a)}V_{k }X  _{ka}\;\Bigg| \;  \max_k|V_{k }|\prec 1, \;  V_{k}\in \cal F_0  \right\}
  \ee
   \end{definition}

\bigskip

 Note: For fixed $k\ne a$, the total number of $X_{ak}$ and $X_{ka}$ ($1\le k\le N$),  in each monomial of the element in $\cal F$ is always even. On the other hand, this number in  $\cal F_{1/2}\cdot \cal F$ is always odd. By the definition, it is easy to see that 
 $$
 \cal F_0, \; \cal F_1\in \cal F 
 $$
$$
\cal F_0\cdot \cal F_{\alpha}=\cal F_{\alpha}, \quad \alpha=0,\;1/2,\;1,\; \emptyset
$$
and 
\be\label{gwyb}
\cal F_{1/2}\cdot\cal F_{1/2}\subset \cal F_{1}, \quad \cal F\cdot\cal F\subset \cal F
\ee

\bigskip

{\bf Examples}: by definition, $G^{(a,a)}_{kl}\le \eta^{-1}$ for any $k,l\neq a$. Hence we have 
$$
\sum_{kl}^{(a)}X _{ka}G^{(a,a)}_{kl}X_{la} \in \eta^{-1} \cal F_1, 
$$ 
$$
\left(\sum_{kl}^{(a)}X _{ka}G^{(a,a)}_{kl}X_{la} \right)\left(\sum_{kl}^{(a)}X _{ak}G^{(a,a)}_{kl}X_{al}\right) \in \eta^{-2} \cal F, 
$$
and  if $\eta=O(1)$
$$
\left(\sum_{kl}^{(a)}X _{ka}G^{(a,a)}_{kl}X_{la} \right)\left(\sum_{kl}^{(a)}X _{ak}G^{(a,a)}_{kl}X_{al}\right)+( \tr X- X_{aa} )\in_n \eta^{-2} \cal F
$$

\bigskip

  \begin{definition}{\bf Uniformness}
 Let $F_T$, $T\in \cal T_N$ be a family of random variables, where $\cal T_N$ is parameter set which may depends on $N$. We say
$$
  F_T\in_n \cal F , \quad T\in \cal T_N, 
  $$
 are {\bf uniform } for all $T\in  \cal T_N$, if the following two uniform conditions hold. 
  \begin{enumerate}
\item There exist uniform integers $m$ and $n$ independent of $N$ such that for all $T\in \cal T_N$, we can write $F_T$ as the sum of $m$ elements in $(\cal F_0\cup \left(\cal F_1 \right)^n)$, i.e., 
  $$
  F_T=\sum_{i=1}^m F_{T, i}, \quad F_{T,i}\in \cal F_0\cup \left(\cal F_1 \right)^n. 
  $$
  \item  All of the stochastic domination relations, i.e., $\prec$,  appearing in all $F_T$'s ($\,T\in \cal T_N$) hold   uniformly. 
 
\end{enumerate}
  
  \medskip
  
    Similarly, for $\cal F_0$, $\cal F_{1/2}$ and $\cal F_{1}$,  we call   $$F_T\in_n \cal F_{\alpha} , \quad T\in \cal T_N, \quad \alpha=0,\;\frac12,\;1  $$ 
   uniformly for all $T\in \cal T_N$, if there exist uniform $m$ independent of $N$ such that $$ F_T=\sum_{i=1}^m F_{T, i}, \quad F_{T, i}\in \cal F_\alpha, \quad \alpha=0,\;\frac12,\;1  $$ and the above  uniform condition  (ii) holds. 
  
  \bigskip
  
  More general,  if $\cal F_\alpha$ is  one of  $\cal F_0,\;\cal F_{1/2},\;\cal F_1,\;\cal F$, so as $\cal F_\beta$,  i.e., $\alpha, \beta =0$, $1/2$, $1$ or $\emptyset$, we say
   $$F_T\in_n \cal F_{\alpha} \cdot \cal F_{\beta},  
   $$ 
     uniformly for all $T\in \cal T_N$ if there exists uniform $m$ independent of $N$ such that $F_T$ can be written as the sum of the $m$ terms in $\cal F_\alpha\cdot \cal F_\beta$, i.e., 
     $$F_T=\sum_{i=1}^mF_{T,\alpha,i}F_{T,\beta,i} $$
     and
     $$
     F_{T,\alpha,i}\in   \cal F_{\alpha}, \quad F_{T,\beta,i}\in  \cal F_{\beta}
     $$
  hold uniformly for all $T\in \cal T_N$.

  \medskip
  
   Furthermore, with fixed $D>0$ and random (or deterministic)  variable  $a_T$, we say 
   $$F_T\in_n  a_T \cal F_{\alpha}\cdot \cal F_{\beta}+O_\prec (N^{-D}),  \quad  \quad  \quad  \quad \cal F_{\alpha},\;  \cal F_{\beta}
    =\cal F_0,\;\cal F_{1/2},\;\cal F_1,\;\cal F $$ 
     uniformly for all $T\in \cal T_N$ if  $F_T$ can be written as 
     $$
     F_T= a_T F_{T,1}+ F_{T,2}
     $$
where  
   $$  F_{T,1}\in _n\cal F_{\alpha} \cdot\cal F_{\beta}, {\quad \rm and\quad  }  F_{T,2}\prec N^{-D} $$
    hold uniformly for all $T\in \cal T_N$.

 \end{definition}

\bigskip
 
Now we estimate the expectation values of   the elements in $\cal F\cdot \cal F_{1/2}$. Let $F_{1/2}\in \cal F_{1/2}$, $F\in \cal F$. With large deviation theory,   we can only obtain 
 \be\nonumber
F_{1/2}\prec 1, \quad   F \prec 1, \quad   F_{1/2}\cdot  F \prec 1
  \ee
But we will show that the elements  in $ \cal F_{1/2}\cdot  \cal F $ may have  much smaller expectation value.
  \begin{lemma}\label{exp22} For fixed indeces $a,b$ and ensemble $X$ in lemma \ref{mainnewL}, 
  let $F_0$ and $F$ be two random variables bounded by $N^C$ for some $C$, i.e., 
  $$
  |F_0|+|F|\leq N^C
  $$ 
  We assume that 
  $$
  F_0\in N^C\cal F_0
 , \quad {\rm   and}\quad  
   F \in_n \cal F_{1/2} \cdot \cal F 
   $$
 Then we have 
 \be\label{EFF1a}
  \left|\E \,F_0 F \right| \prec N^{-1/2}\E\left|F_0\right| +N^{-D} \ee
for any fixed $D>0$. 
  \end{lemma}
  
  
  \bigskip
  
{\it Proof of Lemma \ref{exp22}: } 
For simplicity, we assume $F\in \cal F_{1/2}\cdot \cal F$ (not $\in_n$).  The general case can be proved with the same method. Furthermore, by definition, $ \E F_0F=0$ if $ F\in \cal F_{1/2}\cdot \cal F_0$. 
  Hence one only needs to prove the following case:   for some fixed $m$, $F\in \cal F_{1/2}\cdot (\cal F_1)^m$, i.e., $F$ can be written as the product of one element of $\cal F_{1/2}$ and $m$ elements of $\cal F_1$, i.e., 
  \be\nonumber
  F=\,F_{1/2} F_1F_2F_3\cdots F_m, \quad  F_{1/2}\in \cal F_{1/2} ,\; F_i \in \cal F_1, \quad   1\le i\le m
  \ee
 
By definition, $ F_{1/2} F_1F_2F_3\cdots F_m$ can be consider as a polynomials of $X_{ak}$'s and $X_{ka}$'s ($1\leq k\leq N$), whose coefficients are independents of the $a$-th row and column of $X$.   Then, we can decompose $F$ as 
 \begin{align}\label{412new}
  F= \,&F_{1/2} F_1F_2F_3\cdots F_m 
   \\\nonumber
   =&  
   \sum _{n\leq 2m+1}\;\sum_{k_1,k_2,\ldots, k_n}\;\sum_{s_1,\ldots,s_n}\;\sum_{t_1,\ldots, t_n}\;
 \cal A\Big( \{k_i\}_{  i=1}^n, \{s_i\}_{  i=1}^n,\{t_i\}_{  i=1}^n\Big) \left(\prod_{i=1}^n 
(X_{ak_i})^{s_i}
  (X_{k_ia})^{t_i} \right) 
\end{align}
where $k_i$'s are all different in the summation,  and  $\cal A\Big( \{k_i\}_{  i=1}^n, \{s_i\}_{  i=1}^n,\{t_i\}_{  i=1}^n\Big)$ is the coefficient of  $ \prod_{i=1}^n 
(X_{ak_i})^{s_i}
  (X_{k_ia})^{t_i} $ and it is  independent of the $a$-th row and column of $X$.  We separate the parameter region into two cases.    
 
 {\bf First case: }  $k_i\neq a $ for all $1\le i\leq n$.   By definition of $\cal F_1$, we have 
\be\label{nldd}
\cal A\Big( \{k_i\}_{  i=1}^n, \{s_i\}_{  i=1}^n,\{t_i\}_{  i=1}^n\Big)\prec {\bf 1}\left(\sum s_i+\sum t_i=2m+1\right)\prod_{i=1}^n (N^{1/2})^{\min\{s_i,t_i\}} 
\ee
where the last  factor  come from the $N^{1/2}$   factor in the definition of $\cal F_1$ (see the $N^{1/2}\sum_{k }^{(a)}X_{ak}V_{kk}X_{ka}$ term in the definition of $\cal F_1$. 

{\bf Second case: }  $k_j=a $ for some $1\le j\leq n$. Since the $k_i$'s are all different, hence the other $k_i$'s are not equal to $a$. Let $s_j=s$, $t_j=0$, we have 
\be\label{nldd2}
\cal A\Big( \{k_i\}_{  i=1}^n, \{s_i\}_{  i=1}^n,\{t_i\}_{  i=1}^n\Big) \prec{\bf 1}\left(\sum_{i: \,i\neq j} (s_i+ t_i)\in 2\N+1\right) \prod_{i: \,i\neq j}  (N^{1/2})^{\min\{s_i,t_i\}}N^{s/2} 
\ee

 By definition of $\cal F_1$ and $\cal F$, we know that for any $\delta >0$ and $D>0$, 
 there exists probability set $\Omega$, which is  independent of the $a$-th row and column of $X$,
 such that   $\P (\Omega)\ge 1-N^{-D}$,  and  the $\prec$'s in  \eqref{nldd} and \eqref{nldd2} can be replaced with $\le$.  More precisely, 
\begin{align}\label{nldd3}
 {\bf 1}_\Omega |\cal A_{\rm first\; case}|\leq N^\delta\cdot {\rm  r.h.s} \; {\rm of} \; \eqref{nldd},\quad 
 {\bf 1}_\Omega |\cal A_{\rm second\; case}|\leq N^\delta\cdot  {\rm  r.h.s} \; {\rm of} \; \eqref{nldd2}
 \end{align}  
 With this $\Omega$ and  $|F_0|+|F|\leq N^C$, we have 
\be\label{nldd4}
 \E F_0F= \E  {\bf 1}_\Omega F_0F+ \E  {\bf 1}_{\Omega^c} F_0F=\E  {\bf 1}_\Omega F_0F+O(N^{3C-D})
 \ee
Hence to prove \eqref{EFF1a}, we only need to bound $\E  {\bf 1}_\Omega F_0F $. 
For the first case, i.e., $k_i\neq a$ ($1\leq i\le n$),  using \eqref{nldd3}, and the fact that $F_0$ and $\Omega$ are independent of the $a$-th row and column of $X$,  we have
\begin{align} \nonumber
&\E \sum _{n }\sum_{k_1,k_2,\ldots, k_n }^{(a)}\sum_{s_1,\ldots,s_n}\sum_{t_1,\ldots, t_n}
 {\bf 1}_\Omega F_0  \cal A\Big( \{k_i\}_{  i=1}^n, \{s_i\}_{  i=1}^n,\{t_i\}_{  i=1}^n\Big)
  \left(\prod_{i=1}^n 
(X_{ak_i})^{s_i}
  (X_{k_ia})^{t_i} \right) 
  \\\nonumber
  =&\sum _{n }\sum_{k_1,k_2,\ldots, k_n }^{(a)}\sum_{s_1,\ldots,s_n}\sum_{t_1,\ldots, t_n}
 \E    {\bf 1}_\Omega F_0 \,\cal A\Big( \{k_i\}_{  i=1}^n, \{s_i\}_{  i=1}^n,\{t_i\}_{  i=1}^n\Big)
 \E   \left(\prod_{i=1}^n 
(X_{ak_i})^{s_i}
  (X_{k_ia})^{t_i} \right) \\\nonumber
 \leq 
  & 
   \sum _{n }  \sum_{\{s_i\}}\sum_{\{t_i\}} 
 {\bf 1} \left (\sum s_i+\sum t_i=2m+1\right)\left(\prod_{i=1}^n{\bf 1}\left(  s_i\neq 1\right){\bf 1}\left(  t_i\neq 1\right){\bf 1}\left(  s_i+t_i\neq 0\right)(N^{-1/2})^{\max \{s_i, t_i\}-2}\right) \left(\E|F_0|\right) N^\delta
  \end{align}
 for any $\delta>0$,  where the factor $(N^{-1/2})^{-2}=N^1$ comes from summation of $k_i: 1\leq k_i\leq N$. It is easy to check:   
  \be\nonumber
\prod_i {\bf 1}\left(  s_i\neq 1\right){\bf 1}\left(  t_i\neq 1\right){\bf 1}\left(  s_i+t_i\neq 0\right)(N^{-1/2})^{\max \{s_i, t_i\}-2}\le (N^{-1/2})^{{\bf 1}(s_i+t_i\in 2\N-1)}
  \ee
  Therefore,    for any $\delta>0$, 
\begin{align}\label{kkA1}
  \E& \sum _{n }\sum_{k_1,k_2,\ldots, k_n }^{(a)}\sum_{s_1,\ldots,s_n}\sum_{t_1,\ldots, t_n}
 {\bf 1}_\Omega F_0 \cal A\Big( \{k_i\}_{  i=1}^n, \{s_i\}_{  i=1}^n,\{t_i\}_{  i=1}^n\Big)
  \left(\prod_{i=1}^n 
(X_{ak_i})^{s_i}
  (X_{k_ia})^{t_i} \right)\le \left(\E|F_0|\right) N^{-1/2+\delta}
\end{align}
 Similarly for the second case:  without  loss of generality, we assume $k_1=a$. Then as above, using \eqref{nldd3}, and the fact $\Omega$ independent of the $a$-th row and column of $X$,  we have
\begin{align}\nonumber
&\E \sum _{n }\sum_{ k_2,\ldots, k_n\neq a}\sum_s\sum_{s_2,\ldots,s_n}\sum_{t_2,\ldots, t_n}
 {\bf 1}_\Omega F_0 \cal A\Big( \{k_i\}_{  i=1}^n, \{s_i\}_{  i=1}^n,\{t_i\}_{  i=1}^n\Big)\left(\prod_{i\ne 1}^n 
(X_{ak_i})^{s_i}
  (X_{k_ia})^{t_i} \right) (X_{aa})^s
  \\\nonumber
\le 
  & 
   \sum _{n }  \sum_{\{s_i\}}\sum_{\{t_i\}} {\bf 1}\left(\sum_{i\ge 2} (s_i+ t_i)\in 2\N+1\right) 
   \left( \prod_{i\ge 2}{\bf 1}\left(  s_i\neq 1\right){\bf 1}\left(  t_i\neq 1\right){\bf 1}\left(  s_i+t_i\neq 0\right)(N^{-1/2})^{\max \{s_i, t_i\}+2}  \right)\left(\E|F_0|\right) N^\delta
   \\\label{kkA2}
  \le 
  &  \left(\E|F_0|\right) N^{-1/2+\delta}
\end{align}
Combining \eqref{kkA1} and \eqref{kkA2}, we obtain 
\be\label{jbja}
\E  {\bf 1}_\Omega F_0F\prec \left(\E|F_0|\right) N^{-1/2},
\ee Then together with \eqref{nldd4},  we obtain  \eqref{EFF1a} and complete the proof of Lemma \ref{exp22}. 
  \qed

\bigskip

Now we slightly extend the above lemma. Instead of assuming $  F \in_n \cal F_{1/2}\cdot \cal F$, we assume that $F=  F \in_n \cal F_{1/2}\cdot \cal F+O_\prec (N^{-D})$ for some fixed $D>0$. 

\begin{corollary}\label{exp22+} For fixed indeces $a,b$ and ensemble $X$ in lemma \ref{mainnewL}, 
  let $F_0$ and $F$ be two random variables bounded by $N^C$ for some $C$, i.e., 
  $$
  |F_0|+|F|\leq N^C
  $$ 
  We assume that 
  $$
  F_0\in N^C\cal F_0
$$
  and for some fixed $D>0$, 
  $$
   F =\in_n \cal F_{1/2} \cdot \cal F +O_\prec (N^{-D})
   $$
 Then we have 
 \be\label{EFF1a+}
  \left|\E \,F_0 F \right| \prec N^{-1/2}\E\left|F_0\right| +N^{-D+2C+1} \ee
 \end{corollary}
 
 {\it Proof of Corollary \ref{exp22+}: }   Write 
 $$
 F=F^M+F^e, \quad F^M \in_n \cal F_{1/2} \cdot \cal F, \quad F^e =O_\prec (N^{-D})
 $$
Here superscription $M$ and $e$ are for $main$ and $error$.  (Note $F^M$ and $F^e$ are not assumed to be bounded by $N^C$, otherwise the proof is much simpler.)
 For simplicity, we assume $F^M\in \cal F_{1/2} \cdot \cal F$ (not $\in_n$) and for some $m\ge 0$,  $F^M\in F_{1/2}(F_1)^m$.  Then we repeat the same argument as above, i.e., from \eqref{412new} to \eqref{nldd3}.  
%
 Then  for any (small) $\delta >0$ and (large) $\wt D>0$, 
 there exists probability set $\Omega$, which is  independent of the $a$-th row and column of $X$,
 such that   $\P (\Omega)\ge 1-N^{-\wt D}$,  and \eqref{nldd3} holds. Next we write 
 \begin{align}\nonumber
 \left|\E \,F_0 F \right| &=\left|\E \,1_{\Omega^c} F_0 F \right| +\left|\E \,1_{\Omega} F_0 F^M \right| +\left|\E \,1_{\Omega} F_0 F^e \right|
 \\\nonumber
 &=N^{-\wt D+2C}+N^{-1/2+\delta}\E\left|F_0\right| +\left|\E \,1_{\Omega} F_0 F^e \right|
 \end{align}
 where we used $|F_0|+|F|\leq N^C$  and \eqref{jbja}.   
  
  Now we bound $\left|\E \,1_{\Omega} F_0 F^e \right|$. By the definition of $\prec $ again, there exists $\wt \Omega$ such that $\P(\wt \Omega)\ge 1-N^{-\wt D}$ and 
$$
F^e\leq N^{-D+\delta}
$$
With this $\wt \Omega$, and $|F_0|+|F|\leq N^C$ we write 
\begin{align}\label{1942k}
 \left|\E \,1_{\Omega} F_0 F^e \right|  
 &\leq \left|\E \,1_{\Omega\cap \wt \Omega } F_0 F^e \right|+
\left|\E \,1_{\Omega\cap\wt \Omega ^c } F_0 F^e \right|
\\\nonumber
 &=\left|\E \,1_{\Omega\cap \wt \Omega } F_0 F^e \right|+
\left|\E \,1_{\Omega\cap\wt \Omega ^c } F_0 F  \right|+
\left|\E \,1_{\Omega\cap\wt \Omega ^c } F_0 F^M \right|
\\\nonumber
&\leq 
N^{-D+C+\delta}+N^{-\wt D+2C}+
\left|\E \,1_{\Omega\cap\wt \Omega ^c } F_0 F^M \right|
\end{align}
For the last term, we note that by  the definition of $\Omega$ we can simply bound the term in $\cal F$ which are independent  of the $a$-th row and column  of $X$ by $N^{0.1}$. Then using the assumption $F^M\in \cal F_{1/2}(\cal F_1)^m$,  we have 
$$
|1_{\Omega}F^M|\leq N^{4m+1} \sum^{2m+1} _{n=1} \left(\sum _{k_1, k_2, \cdots, k_{n} }\prod_{j=1}^{2m+1} 
\left(|X _{ak_j}| +|X _{ak_j}|\right) \right)
$$
Together with Cauchy-Schwarz inequality, and subexponential decay property \eqref{subexp},  we obtain that 
$$\left|\E \,1_{\Omega\cap\wt \Omega ^c } F_0 F^e\right|\le 
\E \,1_{\Omega  } \left|F_0 F^e\right|^2\mathbb P(\wt \Omega^c)\leq N^{-\wt D+C_m}
$$
Inserting it into \eqref{1942k}, choosing large enough $\wt D$, we obtain   \eqref{EFF1a+} and complete the proof. 
\qed

\bigskip

More general, if $F_T\in_n \cal F_{1/2}\cdot \cal F$ hold uniformly for $T\in \cal T$, corollary \ref{exp22+} can be extended to the  following integration version. 
\begin{lemma}\label{exp223}
For fixed indeces $a,b$ and ensemble $X$ in lemma \ref{mainnewL}, 
let $F_T$ be a family of random variables such that for some deterministic $x_T$ and uniform $D>0$
  $$
  F_T \in_n  x_T\cal  F_{1/2} \cdot \cal F +O_\prec(N^{-D})
$$
 hold uniformly for $T\in \cal T=\cal T_N$, i.e., $F_{T}=F^M_T+F_T^e$ and 
 $$
  F^M_T \in_n  x_T\cal  F_{1/2} \cdot \cal F, \quad F_T^e= O_\prec(N^{-D})
$$
 hold uniformly for $T\in \cal T=\cal T_N$. 
Here   we assume that $\cup_N\cal T_N$ can be  covered by a compact set in   $ \R^p$ for some $p\in \mathbb N$, this compact set  and $p$ are independent of $N$.

   We also assume that  $|x_T|+|F_T|\leq N^C$ for some uniform $C>0$. Let $F_0$ be  a random variable satisfying $F_0\prec N^C\cal F$ and $|F_0|\leq N^C$.  Then 
 \be\label{EIntFF1a}
  \left|\E \,F_0\int_{T\in \cal T}\, F_T \rd T   \right| \prec N^{-1/2}\left(\E\left|F_0\right|\right) \int_{T\in \cal T}|x_T|  \rd T  +N^{-D+2C+1} \ee 
\end{lemma}
  
 {\it Proof of Lemma \ref{exp223}}: Since $F_0$ and $F_T$ are bounded by $N^C$, one can exchange the order of integration and expectation, i.e., 
 $$
 \E \,F_0\int_{T\in \cal T}\, F_T \rd T=\int_{T\in \cal T} \left(\E \,F_0 \cdot F_T\right) \rd T
 $$
 Then with the uniformness, one can easily extend the proof of Lemma \ref{exp22} and corollary \ref{exp22+}, and prove this lemma. 
 
\qed

 \subsection{Proof of Lemma \ref{mainnewL}. }

The Lem. \ref{exp22} and \ref{exp223} are the key observations for the proof of Lemma \ref{mainnewL}.  Now to  prove Lemma \ref{mainnewL}, we claim that the following lemma \ref{AYFF}, which shows that the terms in  Lemma \ref{mainnewL} 
can be represented by $\cal F$ and $\cal F_{1/2}\cdot \cal F$(with negligible error term).  We first introduce a cutoff function on $\re m^{(a,a)}$. (Recall the definition in Def. \ref{definition of minor})

 \begin{definition}
  Define $\chi_a$ as
\be\label{defchaa}
\chi_a:=\chi_a(\e, w,z)= {\bf 1}\left( |\re m^{(a,a)}|\ge \frac12N^{\e}(N\eta)^{-1}\right)
\ee 
\end{definition}
Note: By definition and \eqref{boundmGm}, $h(t_X)>0$ implies $\chi_a=1$,  and  for any $|\mU|+|\mT|=O(1)$, we have 
\be\label{relhc}
h(t_X)>0 \quad \implies\quad  \chi_a=1\quad \implies \quad |\re m^{(\mU, \mT)}|\ge \frac{1}{4}N^{\e }(N\eta)^{-1}
\ee

\begin{lemma}\label{AYFF} 
Recall $X^{(a,a)}$ and $m^{(a,a)}$ defined in Definition \eqref{definition of minor}. Under the assumption of  Lemma \ref{mainnewL}, for any fixed large $D>0$, we have 
\begin{align}\label{ygfld}
h(t_X)\re m-h(t_{X^{(a,a)}})\re m^{(a,a)}&\in_n \frac{1}{N\eta}\cal F+O_\prec(N^{-D})
\\\label{ygfld2}
   B_m(X) &\in_n  (N\eta)^{m-1}  \cal F+O_\prec(N^{-D}), \quad m=1,2,3
 \\\label{ygfld3}
 \chi_a P_m(X) & \in_n  \frac{1}{N\eta}\cal F_{1/2}\cdot \cal F+O_\prec(N^{-D}), \quad m=1 ,3
 \\\nonumber
\chi_a   P_2(X) & \in_n  \frac{1}{N\eta} \cal F+O_\prec(N^{-D})
\end{align}
 uniformly hold for 
 $$a, b:\;1\leq a\neq b\le N, \quad z: ||z|-1|\le 2\e, \quad {\rm and} \quad w\in I_\e.$$
    \end{lemma}
  
  We postpone the proof of this lemma to the next section. In the remainder of this section, we will prove Lemma \ref{mainnewL} with Lemma \ref{AYFF}. First we introduce a simple lemma for the calculation of  $\cal F$ sets. 

\begin{lemma}\label{lem: prodlaw} Let $A$ and $B$ be two variables stochastically dominated by  $N^C$ 
 for some $C>0$, i.e., $|A|+|B|\prec N^C$. If for random variable $A_0$ and $B_0$, we have 
 $$
 A=A_0+O_\prec(N^{-D}), \quad  B=B_0+O_\prec(N^{-D}), 
 $$
for  some $D>0$. Then
\be\label{2tdgz}
AB=A_0B_0+O_\prec(N^{C-D})
\ee
 \end{lemma}
{\it Proof: } By assumption, $$
\left(A-O_\prec(N^{-D})\right) \left(B-O_\prec(N^{-D})\right)
=A_0B_0
$$
With $|A|+|B|\prec N^C$, we obtain \eqref{2tdgz}.  

\qed
  
  \bigskip

Now we return to finish the proof of Lemma \ref{mainnewL}.

 {\it Proof of Lemma \ref{mainnewL}:  } For simplicity, we introduce the notation $\wt A(w,z)$ as
$$
\wt A(w,z) :=h(t_X)\re m(w,z)-h(t_{X^{(a,a)}})\re m^{(a,a)}(w,z)
$$
First as in \eqref{thpbou}, \eqref{Bbound} and \eqref{PPPNC},  one can see that there exists uniform $C>0$,  such that
\be\label{nnz}
|A^{(f)}_{X }|+|A^{(f)}_{X^{(a,a)} }|+|\wt A(w,z)|+\sum_{n=1,2,3} |\htP_n(X)|+
\sum_{n=1,2,3} |P_n(w,z)|+\sum_{n=1,2,3} |B_n(w,z)|\leq N^C
\ee
With  $A^{(f)}_{X }= A^{(f)}_{X ^{(a,a)}}+ (A^{(f)}_{X }-A^{(f)}_{X ^{(a,a)}})$, we write  $$(A^{(f)}_{X })^{p-3}\htP_1^3(X)=\sum_lC_l(A^{(f)}_{X^{(a,a)} })^{p-3-l}\left( A^{(f)}_{X }-A^{(f)}_{X^{(a,a)} } \right)^l\htP_1^3(X)$$
Recall the definitions in \eqref{defAXf}, \eqref{YYY} and \eqref{defBnn}, for fixed $l$,  with  the notation $\wt A(w,z)$ and \eqref{relhc}, we can write:
\begin{align}\label{jtazs}
&(A^{(f)}_{X^{(a,a)} })^{p-3-l}\left( A^{(f)}_{X }-A^{(f)}_{X^{(a,a)} } \right)^l\htP_1^3(X)
\\\nonumber
=&(A^{(f)}_{X^{(a,a)} })^{p-3-l}N^{l+3}\int_\cal T \prod_{i=1}^l   \wt A(w_i, z_i)\prod_{i=l+1}^{l+3}(\chi_aP_1B_1)(w_i, z_i)
\prod_{i=1}^{l+3}\Delta f(\xi_i)\chi(\eta_i)\phi'(E_i)\rd T
\end{align}
where $\rd T=\prod_i \rd E_i\rd\eta_i\rd A(\xi_i)$ and $\cal T=(I_\e\times \supp f)^{l+3}$. Using \eqref{nnz}, Lemma \ref{AYFF} and Lemma \ref{lem: prodlaw}, for any fixed $D>0$, we have
\be\label{jtazs2}
\prod_{i=1}^l   \wt A(w_i, z_i)\prod_{i=l+1}^{l+3}(\chi_aP_1B_1)(w_i, z_i)\in_n \left(\prod_{i= 1}^{l+3} \frac1{N\eta_i}\right) \cal F_{1/2}\cdot \cal F+O_\prec (N^{-D}), \quad \eta_i=\im w_i
\ee
uniformly hold for $T\in \cal T$. Applying  Lemma \ref{exp223} by choosing 
$$F_0=(A^{(f)}_{X^{(a,a)} })^{p-3-l}, 
\quad 
F_T=\prod_{i=1}^l   \wt A(w_i, z_i)\prod_{i=l+1}^{l+3}(\chi_aP_1B_1)(w_i, z_i),\quad
x_T=\prod_{i= 1}^{l+3}  \eta_i ^{-1}\Delta f(\xi_i)\chi(\eta_i)\phi'(E_i)$$ and $\cal T=(I_\e\times \supp f)^{l+3}$, with \eqref{jtazs} and \eqref{jtazs}, we obtain 
$$
\mathbb E(A^{(f)}_{X^{(a,a)} })^{p-3-l}\left( A^{(f)}_{X }-A^{(f)}_{X^{(a,a)} } \right)^l\htP_1^3(X)
\prec N^{-1/2}\left(\E\left(|A^{(f)}_{X^{(a,a)}  }|^{p-l-3}\right) \right)   +N^{-D } 
$$
for any fixed $D>0$. Then use Holder inequality, we have 
\begin{align}\nonumber
\left|\E (A^{(f)}_{X^{(a,a)} })^{p-3-l}\left( A^{(f)}_{X }-A^{(f)}_{X^{(a,a)} } \right)^l\htP_1^3(X)\right|
 &\prec N^{-1/2}
\left( \OO_\prec(1)+\E\left(|A^{(f)}_{X^{(a,a)}  }|^p\right)\right)
\end{align}
Similarly,  one can prove 
\be\label{hyce2}
\left|\E(A^{(f)}_{X })^{p-3}\htP_1^3(X)\right| +\left| \E(A^{(f)}_{X })^{p-2}\htP_1(X)\htP_2(X)\right|+\left|\E(A^{(f)}_{X })^{p-1}\htP_3(X)\right|\prec N^{-1/2}
\left( \OO_\prec(1)+\E\left((A^{(f)}_{X^{(a,a)}  })^p\right)\right)\ee

It follows from \eqref{boundmGm} that  $m-m^{(a,a)}=O(N\eta)^{-1}$.  Then it is easy to check that $|A^{(f)}_{X^{(a,a)}  }-A^{(f)}_{X}|\le  C$. Inserting it into \eqref{hyce2}, we complete the proof of   Lemma \ref{mainnewL}.

\qed

\section{Polynomialization of Green's functions} 
As showed in the previous sections, to complete the proof of Theorem \ref{z1}, it only remains to prove Lemma \ref{AYFF}. In this section, we will prove Lemma \ref{AYFF}, i.e., write the terms in \eqref{ygfld} as polynomials in $\cal F$ or $\cal F_{1/2}\cdot \cal F$ (up to negligible error). Since the uniformness can be easily checked, we will only focus on the fixed $a, b, z, w$: 
$$a, b:\;1\leq a\neq b\le N, \quad z: ||z|-1|\le 2\e, \quad {\rm and} \quad w\in I_\e.$$

First we need to write the single matrix elements of $G$'s and $\mG's$ as this type of polynomials. To do so, we start with deriving some bounds on $G$'s  under the condition:
 \be\label{conrem}
 |\re m|\ge \frac14N^\e (N\eta)^{-1} 
\ee
Note: this condition is guaranteed by $\chi_a>0$, $h(t_X)>0$ or $h(t_{X^{(a,a)}})>0$.

\subsection{Preliminary lemmas. } This  subsection summarizes some elementary results
from \cite{BouYauYin2012Bulk} and \cite{BouYauYin2012Edge}. 
 Note that all the inequalities in this  subsection hold uniformly for bounded $z$ and $w$. Furthermore, they hold without the condition \eqref{conrem}.

 Recall the definitions of $Y^{(U,T)}$,  $G^{(U,T)}$, $\mG^{(U,T)}$, $\by_i$ and $\mathrm y_i$ in the definition \ref{definition of minor}.

  \begin{lemma} [Relation between $G$, $G^{(\bb T,\emptyset)}$ and  $G^{( \emptyset, \bb T)}$]  \label{lem: GmG}
 For $i,j \neq k  $ ( $i = j$ is allowed) we have
\be\label{111}
 G_{ij}^{(k,\emptyset)}=G_{ij}-\frac{G_{ik}G_{kj}}{G_{kk}}
,\quad
\mG_{ij}^{(\emptyset,k)}=\mG_{ij}-\frac{\mG_{ik}\mG_{kj}}{\mG_{kk}},
\ee
\be\label{Gik}
 G^{ ( \emptyset,i)}
= G+\frac{(G  {\mathrm y} _i^*) \, ( {\mathrm y} _i  G)}
{1-  {\mathrm y} _i G    {\mathrm y} _i ^*}
,\quad
G
 =G^{ ( \emptyset,i)}-\frac{( G^{ ( \emptyset,i)} {\mathrm y} _i^*)  \,
 (  {\mathrm y} _i  G^{ ( \emptyset, i)})}
 {1+   {\mathrm y} _i  G^{ ( \emptyset,i)} {\mathrm y} _i ^* },
 \ee
and
$$
\mG^{ (i,\emptyset)}
=\mG+\frac{(\mG  \by_i) \, (\by_i^* \mG)}
{1-   \by_i^*   \mG     \by_i  }
,\quad
 \mG
 =\mG^{ (i,\emptyset)}-\frac{(\mG^{ (i,\emptyset)}  \by_i)  \,
 ( { \by_i}^ * \mG^{ (i,\emptyset)})}
 {1+  \by_i^*\mG^{ (i,\emptyset)}   \by_i  }.
$$

   \end{lemma}

 \begin{definition}\label{Zi-def}
 In the following, $\E_X$ means the integration with respect to the random variable $X$.
For any $\bb T\subset \llbracket 1,N\rrbracket$, we introduce the   notations
$$
Z^{(\bb T)}_{i }:=(1-\E_{{\mathrm y}_i})
 {\mathrm y}^{(\bb T)}_i  G^{(\bb T, i)} {\mathrm y}_i^{(\bb T)*}
$$
and
$$
\cal Z^{(\bb T)}_{i }:=(1-\E_{\by_i})
\by_i^{(\bb T) *} \mG^{(i, \bb T)} \by_i^{(\bb T)}.
$$
 Recall by our convention that
$\by_i$ is a $N\times 1 $ column vector and $\mathrm y_i$ is a $1\times N $ row vector.
For simplicity we will write
$$
Z _{i }
=Z^ {(\emptyset)}_{i}, \quad \cal Z _{i }
=\cal Z^ {(\emptyset)}_{i}.
$$
\end{definition}

 \begin{lemma}  [Identities for  $G$, $\mG$, $Z$ and  $\cal Z$]   \label{idm}
 For any $ \T\subset \llbracket 1,N\rrbracket$, we have
\begin{align}\label{110}
 G^{(\emptyset , \bb T)} _{ii}
 & =   - w^{-1}\left[1+  m_\mG^{(i, \bb T)}+   |z|^2 \mG_{ii}^{(i, \bb T)} +\cal Z^{(\bb T)}_{i } \right]^{-1},
\end{align}
\be\label{110b}
 {G_{ij} ^{(\emptyset , \bb T) } }
  =   -wG_{ii}^{(\emptyset , \bb T) } G^{(i,\bb T)}_{jj}
\left(   \by_i^{(\bb T)*}  \mG^{(ij, \bb  T)}  \by_j^{(\bb T)}\right) , \quad i\neq j,
\ee
where, by definition,  $\mG_{ii}^{(i,\bb T)}=0$ if $i\in \bb T$.  Similar results hold for $\mG$:
\be\label{110c}
\left[\mG^{(\bb T, \emptyset)} _{ii} \right]^{-1}
  =   - w\left[1+  m_ G^{(\bb  T,i)}+   |z|^2  G_{ii}^{(\bb T,i)} + Z^{(\bb T)}_{i } \right]
\ee
\be\label{110d}
 {\mG_{ij}^{(\bb T, \emptyset)}}
  =   -w\mG_{ii}^{(\bb T, \emptyset)}\mG^{(\bb T, i)}_{jj}\left( \mathrm y_i^{(\bb T)} G^{( \bb T,ij)}   \mathrm y_j^{(\bb T)*}\right), \quad i\neq j.
\ee
\end{lemma}

\begin{definition}[ $\zeta$-High probability events]\label{def:hp}
 Define
\be\label{phi}
\varphi\;\deq\; (\log N)^{\log\log N}\,.
\ee
Let $\zeta> 0$.
We say that an $N$-dependent event $\Omega$ holds with \emph{$\zeta$-high probability} if there is some constant $C$ such that
$$
\P(\Omega^c) \;\leq\; N^C \exp(-\varphi^\zeta)
$$
for large enough $N$.  Furthermore, we say that $\Omega(u)$ holds with \emph{$\zeta$-high probability} uniformly for $u\in U_N$, if there is some uniform constant $C$ such that
\be\label{5yizs}
\max_{ u\in U_N }\P(\Omega^c(u)) \;\leq\; N^C \exp(-\varphi^\zeta)
\ee
for uniformly large enough $N$. 
\end{definition}

Note: Usually we choose $\zeta$ to be 1. By the definition, if some event $\Omega$ holds with $\zeta$-high probability for some $\zeta>0$, then  $\Omega$ holds with probability larger then $1-N^{-D}$ for any $D>0$. 

\begin{lemma}[Large deviation estimate]\label{lem:bh} Let $X$ be defined as in Theorem \ref{z1}. For any $\zeta>0$, there exists $Q_\zeta>0$ such that  for $\bb T\subset\llbracket 1,N\rrbracket$, $|\bb T| \leq N/2$ the following estimates hold  \hp{\zeta} uniformly for $1\leq i,j\leq N$, $ |w|+ |z| \leq C $:
\begin{align}\label{130}
| Z^{(\bb T)}_{i }|=
\left|(1-\E_{ \mathrm y_i})  \left(\mathrm y_i^{(\bb T)} G^{(\bb T,i )}    \mathrm y_i^{(\bb T)*}\right)  \right|
 \leq  \varphi^{  Q_\zeta/2} \sqrt{\frac{\im m_ G^{(\bb T,i )}+ |z|^2 \im G^{(\bb T,i )}_{ii} }{N\eta}}, \\
|\cal Z^{(\bb T)}_{i }| =
 \left|(1-\E_{ \by_i})   \left(\by_i^{(\bb T)*}  \mG^{(i , \bb T )}  \by_i^{(\bb T)} \right) \right|
\leq  \varphi^{  Q_\zeta/2} \sqrt{\frac{\im m_\mG^{(i,\bb T)}+ |z|^2 \im\mG^{(i, \bb T)}_{ii} }{N\eta}}. \non
 \end{align}
 Furthermore, for $i\neq j$, we have
\begin{align}\label{132}
\left|
(1-\E_{\mathrm y_i\mathrm y_j})
\left( \mathrm y_i^{(\bb T)} G^{(\bb T,ij)}   \mathrm y_j^{(\bb T)*}\right)
 \right|
& \leq
  \varphi^{  Q_\zeta/2}
  \sqrt{\frac{\im m_ G^{(\bb T,ij)}
  +|z|^2\im G^{(\bb T,ij)}_{ii}+|z|^2 \im G^{(\bb T,ij)}_{jj}}{N\eta}},
  \\ 
  \left|
(1-\E_{\by_i\by_j} )
\left(\by_i^{(\bb T)*}   \mG^{(ij,\bb T)}  \by_j^{(\bb T)}\right)
 \right|
& \leq
  \varphi^{  Q_\zeta/2}
  \sqrt{\frac{\im m_\mG^{(ij,\bb T)}
  +|z|^2\im\mG^{(ij,\bb T)}_{ii}+|z|^2 \im\mG^{(ij,\bb T)}_{jj}}{N\eta}}, \label{1321}
 \end{align}
 where
 \be\label{1328}
 \E_{\mathrm y_i\mathrm y_j} \left( \mathrm y_i^{(\bb T)} G^{(\bb T,ij)}   \mathrm y_j^{(\bb T)*}\right)
= |z|^2G^{(\bb T,ij)}_{ij}+\delta_{ij}m_G^{(\bb T,ij)}, \quad
\E_{\by_i\by_j}\left(\by_i^{(\bb T)*}   \mG^{(ij,\bb T)}  \by_j^{(\bb T)}\right)= |z|^2\mG^{(ij,\bb T)}_{ij}
+\delta_{ij}m_{\mG}^{(ij,\bb T )}.
 \ee
 \end{lemma}

\bigskip

\begin{lemma}\label{a priori2} Let $X$ be defined as in Theorem \ref{z1}. 
Suppose  $ |w|+ |z| \leq C. $
   For any $\zeta>0$, there exists $C_\zeta$  such that
if the assumption \be\label{eta}
 \eta \geq \varphi^{C_\zeta} N^{-1}|w|^{1/2}
 \ee holds
then the following estimates hold 
 \be\label{53s}
\max_i|G_{ii}|\leq 2(\log N) |w|^{-1/2},
\ee
\be\label{53.5s}
\max_i   |w | |G_{ii}||\mG^{(i,\emptyset )}_{ii}|\leq  (\log N)^4, \ee
 \be\label{54s}
 \max_{ij}|G_{ij}|\leq C(\log N)^2 |w|^{-1/2},
\ee
\be\label{53snew}
|m|\leq 2(\log N) |w^{-1/2}|
\ee
\hp{\zeta} uniformly for   $ |w|+ |z| \leq C $.  
\end{lemma}

 \subsection{Improved bounds on $G\,$'s. }\label{sec:ImbG}\mbox{}

The next lemma gives the bounds on $G$, $\mG$ and  $m$  under the condition \eqref{conrem}.  Note: with \eqref{boundmGm}, 
it  implies that for any $U$, $T$: $|U|+|T|=O(1)$, 
\be\label{conrem33} |\re m^{(U, T)}|\gg (N\eta)^{-1}.
\ee

Before we give the rigorous proof for the bounds on  $G$, $\mG$, we provide a rough picture on the sizes of these terms under the condition \eqref{conrem}, $w\in I_\e$ and $||z|-1|\leq 2\e$.  We note that the typical size of the $G^{(\mU,\mT)}_{kl}$ heavily relies on whether $k=l$ and whether $k$, $l$ are in $\mU$,  $\mT$.  
 \begin{enumerate}
\item If  $k=l\notin \mU\cup \mT $, the typical size of $G^{(\mU,\mT)}_{kk}(w,z)$ is $m(w,z)=\frac1N\tr G(w,z)$. 
\item If  $k\neq l$, and $k, l\notin \mU\cup \mT $, the typical size of $G^{(\mU,\mT)}_{kl}(w,z)$ is 
$\sqrt{|m|/(N\eta)}$. 
\item If $\{k, l\}\cap \mU\neq \emptyset$, then  $G^{(\mU,\mT)}_{kl}=0$.  This result follows from the definition, and it worth to emphasize: 
\be\label{00kl}
\{k, l\}\cap \mU\neq \emptyset\implies  G^{(\mU,\mT)}_{kl}=\mG^{(\mT,\mU)}_{kl}=0 
\ee
\item If $ k=l\in \mT $, then  the typical size of  $G^{(\mU,\mT)}_{kk}$ is  $|wm|^{-1}$
\item If  $k\neq l$, and $ k \in \mT $ and $l\notin \mT$, then  the typical size of  $G^{(\mU,\mT)}_{kl}$ is  
$(|w^{1/2}m |)^{-1} \sqrt{|m|/(N\eta)}$
\item If  $k\neq l$, and $ k,l \in \mT $  then  the typical size of  $G^{(\mU,\mT)}_{kl}$ is  
${|wm^2 |}^{-1} \sqrt{|m|/(N\eta)}$
\item  With the definition 
of $G^{(\mU, \mT)}$ and $\mG^{(\mT, \mU)}$ in Def. \ref{definition of minor}, one can easily see that 
$\mG^{(\mT, \mU)}_{kl}$ has the same typical size as $G^{(\mU, \mT)}_{kl}$(Here the superscript of $\mG$ is $(\mT, \mU)$ not $(\mU, \mT)$).

\end{enumerate}

We note: The $m$ is bounded by $(\log N)^{C}|w|^{-1/2}$ in \eqref{53snew} (no better bound is obtained in this paper), but we believe that it could be much smaller.

\begin{lemma}\label{newBounds}
Let $X$ be defined as in Theorem \ref{z1}.  Let  $\e $ be small enough positive number,   $||z^2|-1|\le 2\e$  and $w\in I_\e$ (see definition in \eqref{defa1a2}).  If \eqref{conrem}  holds , i.e., $|\re m(w,z)|\ge \frac14N^\e ( N\eta)^{-1}$  in $\Omega=\Omega(\e, w, z)$. Then there exists $\wt \Omega\subset \Omega$, and $C>0$ such that  
$\wt\Omega$ holds in $\Omega$ with 1-high probability uniformly for $z$, $w$:  $||z^2|-1|\le 2\e$  and $w\in I_\e$,  (see definition in \eqref{5yizs}) and 
the following bounds hold in $\wt \Omega$ for  any $1\leq i \neq j\leq N$, (Here $A\sim B$ denotes there exists $C>0$ such that $C^{-1}|B|\le |A|\le C|B|$ )
\begin{align}\label{229}
&|1+m |\ge N^{\frac34\e } (N\eta)^{-1}
\\\label{230}
&|1+m^{(i,i)}|\ge N^{\frac14\e}  \left|\cal Z _i^{(i)}\right|
\\\label{231}
& G^{( \emptyset,i)}_{ii}  =(1+O(N^{-\frac14 \e}))\frac{-1}{ w}\frac{1}{1+m^{(i,i)}}
\\\label{232}
&|1+m| \sim |m|
 \\\label{234}
&|G_{ii}| \leq (\log N)^C|m|
\\\label{235} 
 &|G_{ij}^{(\emptyset,i)}| \leq \frac{\varphi^C}{|w^{1/2}m |}\sqrt{\frac{|m|}{N\eta}}
\\\label{234.3}
 &|G^{(\emptyset, j)}_{ii}| \leq (\log N)^C|m|
\\\label{234.7}
 &|G^{( \emptyset, ij )}_{ii}| \leq  \frac{C}{ |w m|}
  \\\label{236}
 &|G_{ij} | \leq \varphi^C  \sqrt{\frac{|m|}{N\eta}}
 \\\label{237}
&|wG_{ii}|^{-1}  \ge   N^{\frac12 \e } |\cal Z_i|
\\\label{238}
&|m^{(i,i)}| \ge (\log N)^{-1}
\end{align}
Furthermore, with the symmetry and the definition 
of $G^{(\mU, \mT)}$ and $\mG^{(\mT, \mU)}$,  
these bounds also hold under the following exchange 
\be\label{exchangeGG}
G^{(\mU, \mT)} \leftrightarrow \mG^{(\mT, \mU)},\quad
\cal Z\leftrightarrow Z.
\ee

 \end{lemma}
{\it Proof of Lemma \ref{newBounds}: } In the following proof, we only focus on the fixed $z$, $w$, $i$ and $j$, since the uniformness can be easily checked.

 We choose $\zeta=1$. Because $\varphi\ll N^\e$ for any fixed $\e>0$ (see \eqref{phi}) and  in this lemma $w\in I_\e$,   one can easily check that  the assumption  in this lemma implies the conditions of   lemma \ref{a priori2} i.e., 
\be\label{relw}
w\in I_\e\quad \implies \eqref{eta}\; {\rm holds\; for}\; \forall C_\zeta
\ee
 Therefore we can use all of the  results (with $\zeta=1$) of  lemma  \ref{a priori2}  in the following proof.
\medskip

{\bf 1.} We first prove \eqref{229}.   The condition \eqref{conrem} implies that $|\frac1N\sum_{i}\re G_{ii}|\ge  \frac14N^{\e}(N\eta)^{-1}$, then there exists $i:$ $1\leq i\leq N$ such that 
$|G_{ii}|\ge \frac14N^{\e}(N\eta)^{-1} $.  Together with  \eqref{53.5s}, it implies that $|\mG^{(i,\emptyset)}_{ii}|\leq  |w|^{-1} N^{-\frac45\e }N\eta$ with $1$ - high probability in $\Omega$. Inserting it into \eqref{110c} with $\mathbb T=i$, using $G^{(i,i)}_{ii}=0$ from \eqref{00kl}, we have  
\be\label{1z78}
|1+m^{(i,i)}+Z_i^{(i)}|\ge N^{\frac45\e }(N\eta)^{-1}
\ee 
Applying  \eqref{130} to bound $Z_i^{(i)}$ with $\mathbb T=i$, using  Schwarz's inequality and the fact $G^{(i,i)}_{ii}=0$ again, we obtain 
\be\label{asZii}
|Z_i^{(i)}|\leq N^{-\e/20} \im m^{(i,i)}+N^{\e/10}(N\eta)^{-1}
\ee
holds with 1-high probability in $\Omega$. Together with \eqref{1z78},   it implies that with 1-high probability in $\Omega$, 
\be\nonumber
|1+m^{(i,i)} |\ge 2N^{\frac34\e }(N\eta)^{-1}
\ee
 Then replacing $m^{(i,i)}$ with $m$ by \eqref{boundmGm}, we obtain  \eqref{229}.  
 
 \bigskip
 
{\bf 2.} {For \eqref{230}}, first using \eqref{boundmGm} and \eqref{229}, we have that for {\it any} $i:$ $1\leq i\leq N$
\be\label{cgya}
|1+m^{(i,i)}|\ge N^{\frac2 3\e  } (N\eta)^{-1}
\ee
holds with 1-high probability in $\Omega$.   Together with the $\cal Z$ version of  \eqref{asZii}:
$$ 
|\cal Z_i^{(i)}|\leq N^{-\e/4} \im m^{(i,i)}+N^{\e/3}(N\eta)^{-1}
 $$
we obtain \eqref{230}. 

\medskip

{\bf 3}.  {For \eqref{231}}, it follows from  \eqref{110} with $\mathbb T=i$, \eqref{00kl} and \eqref{230}. 

\medskip

 {\bf 4.}  {Now we prove \eqref{232}.} Suppose \eqref{229}, \eqref{231} and \eqref{130} holds in $\Omega_0\subset \Omega$. From our previous results,  $\Omega_0$ holds with 1-high probability in $\Omega$.  Now we prove that  \eqref{232} holds in $\Omega_0$. 
 First we assume that   $|1+m|\leq 3$, clearly otherwise \eqref{232}  holds. 
 Together with \eqref{229}, it implies that $(N\eta)^{-1}\le 3N^{-\frac12\e}$.  Using \eqref{boundmGm} and $|1+m|\leq 3$, we obtain $|1+m^{(i,i)}|\leq 4$ 
 and $|m_G^{\emptyset, i)}|\leq 5$. With \eqref{231}, the bound $|1+m^{(i,i)}|\le 4$ implies 
 $|G^{( \emptyset, i)}_{ii}|\ge |5w|^{-1}$.   The assumption $w\in I_\e$ implies $|w|\leq \e$ (see definition of $I_\e$ in \eqref{defa1a2}). 
 Then applying \eqref{130} on $Z_i$, and using $||z|-1|\le 2\e$ and the bounds we just proved on $(N\eta)^{-1}$, $m_G^{( \emptyset, i)}$ and $G^{(\emptyset, i)}_{ii}$, we obtain   that in $\Omega_0$,
\be\label{wdszd}
|Z_i|\leq N^{\frac{-1}3\e}|G^{( \emptyset, i)}_{ii}| 
\ee
  Together with $|G^{( \emptyset, i)}_{ii}| \ge |5w|^{-1}$ and  the assumption $||z|-1|\le 2\e$ and $|w|\le \e$,  we have 
 \be\label{lwtt}
\left||z|^2 |G^{( \emptyset, i)}_{ii}| +Z_{i } \right|\ge |10w|^{-1}
 \ee
 Now inserting \eqref{lwtt} into the identity \eqref{110c}  with  $\mathbb T=\emptyset$, 
using  $|m_G^{( \emptyset, i)}|\le 5$,  and $|w|\le \e$ again,  we obtain that \be\label{237a}
\mG_{ii}= \frac{1}{-w\left(|z|^2 G^{(  \emptyset, i)}_{ii} +Z_{i }\right)} +\e_i, \quad |\e_i|\leq |60w|\frac{1 }{|w|\left(|z|^2 G^{(  \emptyset, i)}_{ii} +Z_{i }\right) }
\ee
Then together with \eqref{wdszd} and \eqref{231}, in $\Omega_0$,  we have
\be\label{phfm}
\left|\mG_{ii}- |z|^{-2}(1+m^{(i,i)})\right|\leq \left(O(|w|)+o(1)\right)| (1+m^{(i,i)})|
\ee
Combining \eqref{229} and \eqref{boundmGm}, we have 
$$
(1+m^{(i,i)})=(1+o(1))(1+m )
$$
Inserting it into \eqref{phfm}, we have 
\be\label{phfm2}
\left|\mG_{ii}- |z|^{-2}(1+m )\right|\leq \left( O(|w|)+o(1)\right)| (1+m )|, \quad {\rm in}\; \Omega_0
\ee
It is easy to extend this result to the following one: 
\be\label{phfm3}
\max_i\left|\mG_{ii}- |z|^{-2}(1+m )\right|\leq \left( O(|w|)+o(1)\right)| (1+m )|, \quad {\rm in}\; \wt \Omega 
\ee
holds  in a probability set $\wt \Omega \subset \Omega$ such that $\wt \Omega $ holds with $1$-high probability in $\Omega$. Since $m=\frac1N\sum _i\mG_{ii}$, for small enough $\e$, with $|w|\le \e$ and $||z^2|-1|\le 2\e$,  \eqref{phfm3}  implies   that 
$$\frac9{10} |1+m| \le |m|\leq \frac{11}{10}|1+m|, \quad {\rm in}\; \wt \Omega $$
It  completed the proof of \eqref{232}.

\bigskip
 
We note: combining \eqref{boundmGm},  \eqref{conrem}, \eqref{229} and \eqref{232}, we have for any $  |U|,\;  |T|=O(1),$
\be\label{longsim}
m^{(U,T)}\sim m\sim 1+m \sim 1+m^{(U,T)}, \quad |U|,\;  |T|=O(1)
\ee

{\bf 5.} {For \eqref{234}}, it follows from  \eqref{231}(with $\mG^{( i,\emptyset)}_{ii}$ in the l.h.s.), \eqref{longsim} and   \eqref{53.5s}.

 \medskip

{\bf 6.} {For \eqref{235}}, first using \eqref{110b}, \eqref{1321}, \eqref{1328} and  \eqref{00kl}, we  obtain that 
\be\label{egll0}
 |G^{(\emptyset, i)}_{ij}|\leq \varphi^C|w| |G^{(\emptyset, i)}_{ii}||G^{(i,i)}_{jj}|
\sqrt{\frac{\im m_ \mG^{(ij,i)}+ |z|^2 \im \mG^{(ij,i )}_{jj} }{N\eta}} 
\ee
holds with 1-high probability in $\Omega$.  Applying  \eqref{53.5s} on $X^{(i,i)}$ instead of $X$, we obtain that   
\be\label{elgg1}
|w||G^{(i,i)}_{jj}||\mG^{(ij,i )}_{jj} | \leq (\log N)^4
\ee
Recall \eqref{conrem} implies \eqref{conrem33}.  Applying  \eqref{234} on $G^{(i,i)}_{jj}$, we have that 
\be\label{elgg2}
|G^{(i,i)}_{jj}| \leq (\log N)^C|m^{(i,i)}|
\ee
holds with 1-high probability in $\Omega$. Then inserting  \eqref{elgg1}, \eqref{elgg2},  \eqref{231} and \eqref{longsim} into \eqref{egll0},  with \eqref{53snew} we obtain \eqref{235}.
  \medskip

{\bf 7.} {For \eqref{234.3}}, from \eqref{Gik}, we have 
$$
G_{ii}
 =G^{ ( \emptyset,j)}_{ii}-\frac{( G^{ ( \emptyset,j)} {\mathrm y} _j^*)_i  \,
 (  {\mathrm y} _j G^{ ( \emptyset, j)})_i}
 {1+   {\mathrm y} _j  G^{ ( \emptyset,j)} {\mathrm y} _j ^* },
$$
On the other hand, \eqref{110c} and \eqref{1328} show that (similar result can be seen in (6.18) of \cite{BouYauYin2012Bulk})$$
\mG_{jj}=-w^{-1}(1+   {\mathrm y} _j  G^{ ( \emptyset,j)} {\mathrm y} _j ^*)^{-1}
$$
Then 
\be\label{my50}
G_{ii}
 =G^{ ( \emptyset,j)}_{ii}+w\mG_{jj}\left(( G^{ ( \emptyset,j)} X^T)_{ij} -G^{ ( \emptyset,j)}_{ij}z^*\right) \,
 \left( ( XG^{ ( \emptyset,j)}  )_{ji}- G^{ ( \emptyset,j)}_{ji}z\right)  
\ee
Since $X_{jk}$'s $(1\le k\le N)$ are independent of $G^{(\emptyset, j)}$,  using large deviation lemma (e.g. see Lemma 6.7 \cite{BouYauYin2012Bulk} ),  as in (3.44) of \cite{BouYauYin2012Edge}, we have that with 1-high probability, 
\be\label{xty1}
|( XG^{ ( \emptyset,j)}  )_{ji}|+|( G^{ ( \emptyset,j)} X^T)_{ij}|\leq \varphi^C\sqrt{\frac{\im G^{ ( \emptyset,j)}_{ii}}{N\eta}}
\ee
Inserting this bound, \eqref{234}, \eqref{235} and \eqref{longsim} into \eqref{my50}, we have
$$
|G_{ii}-G^{ ( \emptyset,j)}_{ii}|\leq  \varphi^Cw|m| \left(\frac{\im G^{ ( \emptyset,j)}_{ii}}{N\eta}+\frac{1}{w|m|N\eta}\right)  
$$
i.e., 
$$
 G_{ii}=\left(1+O(\frac{|w|m}{N\eta})\right)G^{ ( \emptyset,j)}_{ii}+O(\frac{\varphi^C}{ N\eta})  
 $$
 It implies that 
 $$
 G^{ ( \emptyset,j)}_{ii}=\left(1+O(\frac{|w|m}{N\eta})\right)G_{ii}+O(\frac{\varphi^C}{ N\eta})  
 $$
Then with \eqref{53s} and \eqref{53snew}, it implies 
$$
|G_{ii}-G^{ ( \emptyset,j)}_{ii}|\leq \varphi^C(N\eta)^{-1}
$$
and we obtain \eqref{234.3}. 

\medskip

{\bf 8. } For \eqref{234.7}, using \eqref{110} and \eqref{00kl}, we have
$$
G^{ ( \emptyset,ij)}_{ii}=-w^{-1} [1+m_{\mG}^{(i,ij)}+\cal Z_i^{(ij)} ]^{-1}
$$
 Using \eqref{130} and \eqref{00kl} again, we can bound $\cal Z_i^{(ij)}$ as 
 $$
| \cal Z_i^{(ij)}|\leq \varphi^C\sqrt{\frac{\im m_\mG^{(i,ij)}}{N\eta}}
 $$
 Together with \eqref{longsim} and \eqref{229}, we obtain \eqref{234.7}. 
 
 \medskip
 
 {\bf 9. }{For \eqref{236}},    using \eqref{110b}, \eqref{1321} and \eqref{1328}, we  obtain that 
\be\label{egzz}
 | G _{ij}|\leq \varphi^C|w| | G_{ii}|| G^{(i,\emptyset)}_{jj}|
\sqrt{\frac{\im m_ \mG^{(ij,\emptyset)}+ |z|^2 \im \mG^{(ij ,\emptyset)}_{jj}+ |z|^2 \im \mG^{(ij ,\emptyset)}_{ii}}{N\eta}} +
\varphi^C|wz^2|| G_{ii}|| G^{(i,\emptyset)}_{jj}||\mG^{(ij ,\emptyset)}_{ij}|
\ee
Furthermore, with \eqref{110d}, \eqref{132},  \eqref{00kl} and \eqref{longsim}, we have 
\be\label{egzz2}
|\mG^{(ij ,\emptyset)}_{ij}|\leq \varphi^C|w||\mG^{(ij ,\emptyset)}_{ii}||\mG^{(ij ,i)}_{jj}|
\sqrt{\frac{\im m^{(ij,ij)}   }{N\eta}} 
\leq \varphi^C|w||\mG^{(ij ,\emptyset)}_{ii}||\mG^{(ij ,i)}_{jj}|
\sqrt{\frac{|m|   }{N\eta}} 
\ee
Here these two bounds holds with 1-high probability. As in \eqref{elgg2}, applying  \eqref{231} on $\mG^{(ij,i  )}_{jj}$, with \eqref{longsim} we have 
$$
 |\mG^{(ij,i  )}_{jj}|\leq C|w|^{-1}|m^{(ij,ij)}|^{-1}\leq C|w|^{-1}|m|^{-1}
$$
with 1-high probability in $\Omega$. With \eqref{234}, \eqref{234.7}, \eqref{boundmGm} and \eqref{229}, we also have 
$$
| G_{ii}|\leq(\log N)^C|m| 
,\quad
|\mG^{( ij,\emptyset )}_{ii}|+|\mG^{( ij,\emptyset )}_{jj}|\leq C|w|^{-1}|m |^{-1},\quad  |m_ \mG^{(ij,\emptyset)}| \leq C|m|, 
\quad
 $$
  For the $G_{jj}^{( i,\emptyset)}$ in \eqref{egzz},  as in \eqref{my50} and \eqref{xty1}, with \eqref{00kl}, we have 
  \begin{align}\label{jysr}
G_{jj}^{( i,\emptyset)}
-G_{jj}^{( i,i)}&=w\mG_{ii}^{( i,\emptyset)} ( G^{ ( i,i)} X^T)_{ji}   ( XG^{ ( i,i)}  )_{ij} 
\\\nonumber &=O\left(\varphi^C |w\mG_{ii}^{( i,\emptyset)} | \im G_{jj}^{( i,i)} (N\eta)^{-1}\right)
  &
\end{align}
Then applying \eqref{234} on $G_{jj}^{( i,i)}$, and applying \eqref{231} on $\mG_{ii}^{( i,\emptyset)}$, with \eqref{longsim} we obtain that 
 $$| G_{jj}^{( i,\emptyset)}|\leq(\log N)^C|m | $$
Inserting these bounds   into \eqref{egzz} and \eqref{egzz2}, we obtain \eqref{236}. 
 
 \bigskip
 
 {\bf 10.} {For \eqref{237}}, using \eqref{130} (with $\mathbb T=\emptyset$) and  \eqref{231}, \eqref{longsim} , we have 
 \be\label{cza}
|  \cal Z_i|\leq \varphi^C\sqrt{\frac{ |m|+(|w m|)^{-1}}{N\eta}}
 \ee
 holds with 1-high probability in $\Omega$. Together with   \eqref{53snew}, we obtain 
 \be\label{cza2}
|  \cal Z_i|\leq \varphi^C\sqrt{\frac{  (|w m|)^{-1}}{N\eta}}
 \ee
Together with  \eqref{234} and  \eqref{53snew}, we have  
$$|  \cal Z_i| | wG_{ii}|  \leq \varphi^C\sqrt{\frac{|w|^{1/2}}{N\eta}}.$$ 
Then with \eqref{wzzln},  we obtain \eqref{237}.
 
 \medskip
 
{\bf 11.} {For \eqref{238}, }  we note that  \eqref{232} implies $|m|\ge (\log N)^{-1}$.  Then with \eqref{longsim}, we obtain \eqref{238}. 
 
  \qed

\subsection{Polynomialization of Green's functions: }
 
 In this subsection, using the bounds we proved in the last subsection, we write the $G$'s and $\mG$'s as the polynomials in $\cal F$ and $\cal F_{1/2}\cdot \cal F$ (with negligible error). 
 
 We note: In the  Lemma \ref{mainnewL} and \ref{AYFF} we assumed $X_{ab}=0$, but  the bounds we proved in Lemma \ref{a priori2} and Lemma \ref{newBounds} still hold for this type of $X$, the similar  detailed argument was given in Remark 3.8 of [2].    
 
 \begin{lemma} Lemma \ref{a priori2} and Lemma \ref{newBounds} still hold if one enforces $X_{st}=0$ for some fixed $1\leq s,t\leq N$. 
 \end{lemma}
 Note: Here $s,t$ are allowed to be the same as the $i,j$ in Lemma \ref{a priori2} and Lemma \ref{newBounds}.  For example, from \eqref{236}, we have $|G_{st}|\leq \varphi^Cm^{1/2}(N\eta)^{-1/2}$, even if $X_{st}=0$.
 \bigskip

 By the definitions of $A_X^{(f)}$, $\htP_{1,2, 3}(X)$, $B_{1,2,3}(X)$ and $P_{1,2,3}(X)$, one can see that the values of $A_X^{(f)}$, $\htP_{1,2, 3}(X)$ would not change if one replaced the $G$'s inside with  $\chi_a G$'s.
Therefore, instead of   $G$'s, we will write  $\chi_aG$  as the polynomials in $\cal F$ and $\cal F_{1/2}\cdot \cal F$ (with negligible error).
\begin{definition}For simplicity, we define the notations:
\be\nonumber
\alpha:=\chi_a |m^{(a,a)}|,\quad  \beta:=\frac{\chi_a}{|w m^{(a,a)}|}, \quad \gamma=\chi_a|w|^{1/2}\sqrt{\frac{|m^{(a,a)}|}{N\eta}}
\ee
\end{definition}
We collect some basic properties of these quantities in the the following lemma. 
\begin{lemma}\label{yzb}
Under the assumption of Lemma \ref{mainnewL},  for $z$, $w$:  $||z^2|-1|\le 2\e$  and $w\in I_\e$
\begin{align}\label{abcB}
& \chi_a(\log N)^{-1}\leq \alpha \leq  (\log N)^C\beta\leq (\log N)^C\eta^{-1}\\
\label{abcB2}
 & \chi_a(\log N)^{-1}N^{-1/2}\leq \gamma\leq N^{- \e/ 2}
 \\\label{abcB3}
 & \beta\gamma ^2=\chi_a(N\eta)^{-1} 
 \\\label{abcB4}
 &\frac{\chi_a(\log N)^C}{N\eta}\leq\alpha\leq \chi_a(\log N)^C|w^{-1/2}|
 \end{align}
 hold with 1-high probability.
 \end{lemma}
 
 {\it Proof of Lemma \ref{yzb}: } We note $\chi_a=1$ implies the condition \eqref{conrem}. Hence  the results in  Lemma \ref{newBounds} hold with 1- high probability.  First from \eqref{238} and $|w|\ge \eta$, we have  the first and the third inequalities of \eqref{abcB}, and the first inequality of \eqref{abcB2}. The second inequality in \eqref{abcB} follows from \eqref{53snew} and \eqref{longsim}.  It also implies the second inequality of \eqref{abcB4}.   Combining  the second inequality of \eqref{abcB} with \eqref{wzzln}, we obtain the second inequality in \eqref{abcB2}. For \eqref{abcB3}, one can easily check this identity by the definition of $\beta$ and $\gamma$.  For the first inequality of \eqref{abcB4}, it follows from
 \eqref{229} and \eqref{longsim}.
  \qed
\begin{definition}\label{defSS}  Under the assumption of Lemma \ref{mainnewL},  for $w\in I_\e$, $||z|-1|\leq 2\e$ and $s, k\neq a$, we define $S_{ks}$ and $\wt S_{sk}$ as random variables which are independent of the $a$-th row and columns of $X$ and 
$$\frac{ G ^{( \emptyset, a)}_{ka}}{ G ^{( \emptyset, a)}_{aa}}
=\sum_{s}^{(a)}S_{ks}X_{sa}\quad {\rm and }
\quad \frac{ G ^{( \emptyset, a)}_{ak}}{ G ^{( \emptyset, a)}_{aa}}=\sum_{s}^{(a)}X_{sa}\wt S_{sk}
$$
With  \eqref{110b}, one can obtain their explicit expressions, e.g., 
\be\nonumber
S_{ks}:=z^* wG^{(a,a)}_{kk}\mG^{(ak,a)}_{ks}-wG^{(a,a)}_{kk}\sum_{t}^{(a)}\mG^{(ak,a)}_{st}X_{tk}
\ee
Similarly, we define $\cal S_{ks}$ and $\wt {\cal S}_{sk}$ as random variables which are independent of the $a$-th row and columns of $X$ and 
$$\frac{ \mG ^{( \emptyset, a)}_{ka}}{ \mG ^{( a, \emptyset)}_{aa}}
=\sum_{s}\cal S_{ks}X_{as}\quad {\rm and }
\quad 
\frac{ \mG ^{( \emptyset, a)}_{ak}}{ \mG ^{( a, \emptyset)}_{aa}}
=\sum_{s}X_{as}\wt {\cal S}_{sk}
$$

\end{definition}
 As one can see that $S$, $\wt S$, $\cal S$ and $\wt {\cal S}$ have the same behaviors. Here we collect some basic properties of these quantities in the the following lemma. 
\begin{lemma}\label{jsds}
We assume  that  $||z|-1|\leq 2\e$, $w\in I_\e$, $k\neq a$ and $X$ satisfies the assumption of Lemma \ref{mainnewL}.  For some $C>0$, with 1- high probability, we have 
\be\label{boundS}
 |\chi_aS_{ks}|  \leq \chi_a\varphi^C\left(\delta_{sk}+\gamma \right)
\ee
so as $\wt S$, $\cal S$ and $\wt {\cal S}$. Recall the definition $\cal F$'s in Def. \ref{def: FFF},  for some $C>0$, we have
\be\label{xsx}
\chi_aX_{aa}\in_n    \gamma \cal F,\quad \chi_a(XSX)_{aa}\in_n    \gamma \cal F,
\ee
and
\be\label{xsx2}
 \chi_a(X^T\wt S SX)_{aa}\in_n   N\gamma^2\cal F
\ee
Furthermore, \eqref{boundS}, \eqref{xsx} and \eqref{xsx2} hold uniformly for $||z|-1|\le 2\e$, $w\in I_\e$ and  $k,s:   k,s\neq a$, $1\leq k,s \le N$. 
\end{lemma}
Note:  With \eqref{xsx}, we also have \be\label{xsx3}
\chi_a\left( G_{aa} ^ {( \emptyset, a)}\right)^{-1}( XG ^ {( \emptyset,a)})_{aa}
= \chi_a\left( G_{aa} ^ {( \emptyset, a)}\right)^{-1}\sum _{k} X_{ak}G_{ka}^ {(\emptyset,a)} =\chi_a\left((XSX)_{aa}+X_{aa}\right)  \in_n   \gamma \cal F
 \ee

{\it Proof of Lemma \ref{jsds}: }  
Since the uniformness are easy to be checked, we will only focus on the fixed $z$, $w$, $s$ and $k$.

\medskip
 
{\bf 1.} {For \eqref{boundS}},  the condition $\chi_a=1$ implies that we can apply Lemma \ref{newBounds} on the  $X^{(a,a)}$. Recall: these bounds also hold under the exchange \eqref{exchangeGG}. Then the bounds \eqref{234} and \eqref{235} imply that for $s\neq k$, 
\be\label{xyss-1}
\chi_a| G^{(a,a)}_{kk}| \leq (\log N)^C|m^{(a,a)}|, 
\quad 
\chi_a|\mG^{(ak,a )}_{ks}| \leq \frac{\varphi^C}{|w^{1/2}m^{(a,a)} |}\sqrt{\frac{|m^{(a,a)}|}{N\eta}}, 
\ee
 holds with 1-high probability.  Similarly \eqref{231} and \eqref{longsim} implies that for $s=k$
\be\nonumber
\chi_a|\mG^{(ak,a )}_{kk}| \leq C|w m^{(a,a)} |^{-1}  
\ee
 holds with 1-high probability. Then with the explicit expression of $S_{ks}$ in Def. \ref{defSS}, we have
\be\label{xyss0}
 \chi_aS_{ks}= O(\delta_{ks}+\varphi^C\gamma )-wG^{(a,a)}_{kk}\sum_{t}^{(a)}\mG^{(ak,a)}_{st}X_{tk}
 \ee
 holds with 1-high probability.    Since $X_{tk}$'s are independent of $\mG^{(ak,a)}_{st}$'s ($1\leq t\le N$),  using large deviation lemma (e.g. see Lemma 6.7 \cite{BouYauYin2012Bulk} ),  as in (3.44) of \cite{BouYauYin2012Edge}, we have for 
 \be\label{xyss1}
 |\sum_{t}^{(a)}\mG^{(ak,a)}_{st}X_{tk}|\leq \varphi^C\sqrt{\frac{\im \mG^{(ak,a)}_{ss}}{N\eta}}
 \ee
 holds with 1-high probability.   Applying  Lemma \ref{newBounds}
  on the  $X^{(a ,a )}$ again,  from \eqref{234.3}, we have 
  $$
 | \mG^{(ak,a)}_{ss}| \leq (\log N)^C|m^{(a ,a )}|+C\delta_{ks} \frac{1}{|w m^{(a,a)}|}
 $$
  with 1-high probability.   Together with  the first part of \eqref{xyss-1}, \eqref{xyss0} and  \eqref{xyss1}, we obtain 
 \be\label{xyss3}
|\chi_aS_{ks}|\leq C\delta_{sk}+\varphi^C\gamma +\varphi^C| w^{1/2}m^{(a,a)}|\gamma
\ee
 with 1-high probability.  At last, with \eqref{53snew} and \eqref{longsim}, we obtain \eqref{boundS}.  
 
 \medskip
 
 {\bf 2.} {For \eqref{xsx}}, we recall  the definition of $\cal F$ in Def. \ref{def: FFF}, especially the two $N^{1/2}$ factors in $\cal F$. It is easy to see that \eqref{xsx} follows from   the first inequality of \eqref{abcB2} and the bounds on $S$ in \eqref{boundS}.

\medskip

 {\bf 3.} {For \eqref{xsx2}}, since the \eqref{boundS} also holds for  $\wt S$, then with the  first inequality of \eqref{abcB2}, we have 
 $$
|\chi_a(\wt SS)_{kl}|\le \varphi^C \left(\delta_{kl}+\gamma+N\gamma^2\right)\leq \varphi^CN\gamma^2
 $$
 with 1-high probability. Together with definition of $\cal F$, we obtain \eqref{xsx2}. 

   \qed

 \bigskip

Now we introduce a method to track and show the dependence of the random variables on  
 the indices. First we give a simple example to show the basic idea. Let $A_{kl}$, $1\leq k,l\le N$ be a family of random variables:
 \be\label{ydddgd}A_{kl}= \frac{G^{(a,a)}_{kk}}{|G^{(a,a)}_{kk}|}\frac{G^{(a,a)}_{ll}}{|G^{(a,a)}_{ll}|} (X G^{(a,a)}X^T)_{aa}, \quad 1\leq k,l\le N\ee
where $X^T$ is the transpose of $X$. By definition of $\cal F$ and $\cal F_0$, we can say, 
$$
A_{kl}\in \cal F_0\cdot \cal F_0\cdot \cal F_1 \in \cal F
$$
But  the first part of the r.h.s. of \eqref{ydddgd}, i.e., $\frac{G^{(a,a)}_{kk}}{|G^{(a,a)}_{kk}|}$ only depends on the first index $k$, the second part  $\frac{G^{(a,a)}_{ll}}{|G^{(a,a)}_{ll}|}$ only depends on the second index $l$ and the third part is independent of the indices. Therefore, we prefer to write it as 
$$
A_{kl}\in \cal F_0^{[k]}\cdot \cal F_0^{[l]}\cdot\cal F_1^{[\emptyset]} .
$$
 More precisely, $A_{kl}\in \cal F_0^{[k]}\cdot \cal F_0^{[l]}\cdot\cal F_1^{[\emptyset]} $ means that 
$
A_{kl}=f_1(k)f_2(l)f_3
$ and $f_1(k)\in \cal F_0$, $f_2(l)\in \cal F_0$, $f_3\in \cal F_1$, and $f_{1}(k)$ only depends on index $k$,  $f_{2}(l)$ only depends on index $l$,  and   $f_{3} $ does not depends on index.

For general case, to show how the variable depends on the indices, we define the following notations. 
\begin{definition}\label{def:Fk}  Let $A_I$ be a family of  random 
variables  where $I$ is  indices (vector), not including index $a$.     we write 
$$
A_I\in \prod_i \cal F^{[I_i]}_{\alpha_i}, \quad \cal F _{\alpha_i}\in \left\{\cal F_0, \;\cal F_{1/2}, \;\cal F_1,\;\cal F\right\}
$$ 
where $ I_i$ is a part of $I$, if and only if there exists $ f_{ i}(I_i)\in  \cal F _{\alpha_i}$  such that $A_I=\prod_i f_{ i}(I_i)$ and $f_{ i}(I_i)$ only depends on the indices 
 in $ I_i $. \end{definition}
For the example in \eqref{ydddgd}, we write $A_{kl}\in \cal F_0^{[k]}\cdot \cal F_0^{[l]}\cdot\cal F^{[\emptyset]} $, where 
$I=(k, l)$, $I_1=(k)$, $I_2=(l)$ and $I_3=(\emptyset)$, $\alpha_1=\alpha_2=0$ and $\alpha_3=1$
\bigskip

The following lemma shows the $G$'s can be written as the polynomials in $\cal F$'s. 

\begin{lemma}\label{boundchaa}For simplicity, we introduce the notaion:
\be\label{gniz}
 \cal F_{ 0, X}^{[k]}  :=X_{ak} \cal F^{[k]}_0+ X_{ka }
\cal F^{[k]}_0  
\ee
i.e.,  
$$
f_k\in  \cal F_{ 0, X}^{[k]}\iff \exists g_k , h_k\in \cal F^{[k]}_0:  f_k=X_{ak}g_k  + X_{ka }h_k
$$

 Let $w\in I_\e$ and $||z|-1|\le2\e$.  Under the assumption of Lemma \ref{mainnewL},  for any  $D>0$, we have 

  \be\label{cGaaa}
\chi_a G^{( \emptyset, a)}_{aa}, \in_n  \beta\cal F+O_\prec(N^{-D})
 \ee
  and
  \be\label{cGaa}
\chi_a G_{aa} \in_n   \alpha \cal F   +O_\prec(N^{-D})
\ee
  For any $k\neq a$,   
 \be\label{cGaab}
 \chi_a (G^{(\emptyset, a)}_{aa})^{-1}G^{(\emptyset, a)}_{ka}\in_n
   \gamma  \cal F_{1/2}^{[k]}
 + \cal F_{ 0, X}^{[k]}
 +O_\prec(N^{-D}) 
 \ee
 and,    
\be\label{cGab}
 \chi_a G_{ak} \in_n\sqrt{\frac{\alpha}{N\eta}}  \cal F_{1/2}^{[k]}\cdot \cal F ^{[\emptyset]}
+\left(\alpha + \beta\gamma\right) \cal F^{[\emptyset]}\cdot \cal F_{ 0, X}^{[k]}+O_\prec(N^{-D}) 
 \ee
 For  any
  $k, l\neq a$, 
\begin{align}\label{cGkl}
\chi_a\left(G_{kl}-G^{(a,a)}_{kl}\right)\in _n &\;
\Bigg(\frac{\chi_a}{N\eta} \cal F^{[k]}_{1/2}  \cal F^{[l]}_{1/2} 
+
\beta \gamma  \cal F_{ 0, X}^{[k]} \cal F_{1/2}^{[l]}
+\beta\gamma \cal F_{ 0, X}^{[l]} \cal F_{1/2}^{[k]}
  +\beta \cal F_{ 0, X}^{[k]}  \cal F_{ 0, X}^{[l]}  \Bigg)\cal F^{[\emptyset]}+O_\prec(N^{-D})
 \end{align}
Furthermore, \eqref{cGaaa}-\eqref{cGkl} hold uniformly for $||z|-1|\le 2\e$, $w\in I_\e$ and $1\leq k,l\neq a\le N$ \end{lemma}

{\it Proof of Lemma \ref{boundchaa}: } Because one can easily check the uniformness, in the following proof we will only focus on the fixed $w$, $z$, $k$ and $l$. Recall  \eqref{relhc} and \eqref{relw},  with the assumption $w\in I_\e$ and $||z|-1|\le2\e$, we know the  results in Lemma \ref{a priori2}  and \ref{newBounds} hold under the assumption of this lemma.   Furthermore, these results also  hold  for $X^{(a,a)}$(instead of $X$). 
 
{\bf 1.} {We first prove  \eqref{cGaaa}}.  Applying  Lemma \ref{newBounds} on $X^{(a,a)}$, with  
\eqref{234}, \eqref{236} and the first inequality of \eqref{abcB4}, we have 
\be\label{qbq}
 \chi_aG^{(a,a)}_{kl }\in \left(\delta_{kl}\,\alpha+|w^{-\frac12}|\gamma\right)\cal F_0, \quad {\rm and }\quad |w^{-\frac12}|\gamma \leq \alpha
 \ee Then with  \be\label{fhjd}
 Z_{a}^{(a)}=\left((XG^{(a,a)}X^T)_{aa}-m^{(a,a)} \right)
 \ee and $\alpha:=\chi_am^{(a,a)}$, we have
\be\label{qbq2}
 \chi_a(XG^{(a,a)}X^T)_{aa}\in_n \alpha\, \cal F  \quad{\rm and}\quad  \chi_aZ_{a}^{(a)}\in_n \alpha\,\cal F. 
\ee
From  \eqref{110c} and \eqref{00kl} with $i=a$, $\mT=a$, we have
$$
\chi_a\mG^{(a, \emptyset)}_{aa}=\chi_a\frac{1}{-w}\frac{1}{1+m^{(a,a)}+Z_a^{(a)}}
$$
Then  with   \eqref{230},  for any $\e, D>0$, there exists $C_{\e, D}$ depending $\e$ and $D$, such that 
\be\nonumber
\chi_a\mG^{(a, \emptyset)}_{aa}=\chi_a\frac{1}{-w}\sum_{k=1}^{C_{\e, D}}\left(\frac{1}{(1+m^{(a,a)})^k} (  Z_a^{(a)})^{k-1} \right)+O_\prec(N^{-D})
\ee
holds with 1-high probability. 
 Hence with \eqref{longsim} and  $\chi_aZ_{a}^{(a)} \in_n m^{(a,a)} \cal F$ in \eqref{qbq2}, we obtain that
\be\label{mttaz}
\chi_a\mG^{(a, \emptyset)}_{aa }\in_n \frac{1}{wm^{(a,a)}}\cal F+O_\prec(N^{-D})=\beta  \cal F +O_\prec(N^{-D})
 \ee
 which implies \eqref{cGaaa} with the fact: $\mG^{(a, \emptyset)}_{aa }$ and $ G^{(  \emptyset, a)}_{aa }$ have the same behavior.

 {\bf 2.} { Now we prove \eqref{cGaa}}. From \eqref{237} and \eqref{110}, with $i=a$ and $T=\emptyset$,  for any $\e, D>0$, there exists $C_{\e, D}$ depending $\e$ and $D$ such that with 1-high probability, 
\be\label{tjjs}
 \chi_aG_{aa}=\sum_{k =1}^{C_{ \e, D}}\frac{-w^{-1} \chi_a}{(1+  m_\mG^{(a, \emptyset)}+   |z|^2 \mG^{(a, \emptyset)}_{aa})^k}(\cal Z_a)^{k-1} +O_\prec(N^{-D})
\ee
Note: $1+  m_\mG^{(a, \emptyset)}+   |z|^2 \mG^{(a, \emptyset)}_{aa}$ is independent of the $a$-th column of $X$, but depends on the $a$-th row of $X$. From \eqref{130} and \eqref{1328}, we have 
\be\label{rbhs}
\cal Z_a= z \sum _{k}(X^T)_{ak}\mG^{(a,\emptyset)}_{ka} + z ^*\sum _{k}\mG^{(a,\emptyset)}_{ak}X_{ka}+ 
\sum _{kl}(X^T)_{ak}\mG^{(a,\emptyset)}_{kl}X_{la}- m_\mG^{(a,\emptyset)}-|z|^2 \mG^{(a,\emptyset)}_{aa}
\ee
 Now we claim that for any $D$, 
\be\label{zGF}
\chi_a\cal Z_a \in_n \beta\cal F  +O_{\prec}(N^{-D})
\ee
and 
\begin{align}\label{DDs3}
 \chi_a \left({1+  m_\mG^{(a, \emptyset)}+   |z|^2 \mG^{(a, \emptyset)}_{aa} }\right)^{-1}
  &\in_nw\alpha\,\cal F+O_\prec(N^{-D})
\end{align}
Combining \eqref{zGF}, \eqref{DDs3} and \eqref{tjjs}, we obtain \eqref{cGaa}. 

{\bf 2.a} {We prove \eqref{zGF}} first.  Using the $\mG$ version of \eqref{xsx3} and \eqref{mttaz}, we can write the first two terms of the r.h.s. of \eqref{rbhs} as    we can write
\begin{align}\label{23cl}
\chi_a z \sum _{k}(X^T)_{ak}\mG^{(a,\emptyset)}_{ka} + \chi_a z ^*\sum _{k}\mG^{(a,\emptyset)}_{ak}X_{ka}=
2\chi_a\re z\,\left(\mG ^{(a,\emptyset)} X\right)_{aa}\in_n \beta \gamma \cal F+O_\prec(N^{-D})
\end{align}
Similarly for the third term of the  r.h.s. of \eqref{rbhs}, using \eqref{111},  we can write it as
 \begin{align}\nonumber
( X^T \mG^{(a,\emptyset)}X )_{aa}&
 =\sum _{kl}(X^T)_{ak}\mG^{(a,a)}_{kl}X_{la} 
+ (\mG ^{(a,\emptyset)}_{aa})^{-1} \sum_{kl}(X^T)_{ak} \mG_{ka}^ {(a,\emptyset)} \mG^{(a,\emptyset)} _{al}X_{la} 
\\\nonumber
=&( X^T \mG^{(a,a)}X )_{aa}+ (\mG ^{(a,\emptyset)}_{aa})^{-1} \left(\mG ^{(a,\emptyset)} X\right)_{aa}\left(\mG ^{(a,\emptyset)} X\right)_{aa}
\end{align}
Using \eqref{qbq2} , \eqref{xsx3} and \eqref{mttaz},  we obtain 
 \be\label{24cl}
\chi_a ( X^T \mG^{(a,\emptyset)}X )_{aa}\in _n \alpha \cal F+\beta \gamma ^2\cal F+O_\prec(N^{-D}), \quad 
\ee
For the fourth term of the  r.h.s. of \eqref{rbhs}, using  \eqref{111}, we have 
$$\mG_{kk}^{(a,a)}=\mG_{kk}^{(a,\emptyset)}-\frac{\mG_{ka}^{(a,\emptyset)}\mG_{ak}^{(a,\emptyset)}}{\mG_{aa}^{(a,\emptyset)}}
$$
Together with \eqref{xsx2}, it implies that

\be\label{wxsl}
m_\mG^{(a,\emptyset)}=m^{(a,a)} 
+\frac{1}{N} \mG ^{(a,\emptyset)}_{aa} \left( \left(  X \wt {\cal S}\cal S X^T \right)_{aa}+1\right)
\ee
and 
$$\chi_am_\mG^{(a,\emptyset)}\in _n \alpha \cal F+\beta \gamma ^2\cal F+\frac1N\beta \cal F
+O_\prec(N^{-D})$$
Now inserting these bounds back to \eqref{rbhs}  and using the relations between $\alpha$, $\beta$ and $\gamma$ in \eqref{abcB} and \eqref{abcB2}, we 
  obtain \eqref{zGF}.  

\bigskip

{\bf 2.b} { Now we prove \eqref{DDs3}}. With \eqref{wxsl} and 
$$(\mG ^{(a,\emptyset)}_{aa})^{-1}=-w(1+(XG^{(a,a)}X^T)_{aa})=-w(1+m ^{(a, a)}+Z_{a}^{(a)})
$$ we write
\begin{align}\label{ydnn}
\frac{1}{1+  m_\mG^{(a, \emptyset)}+   |z|^2 \mG^{(a, \emptyset)}_{aa} }
&=\frac{\left(\mG ^{(a,\emptyset)}_{aa}\right)^{-1}}{\frac{1+m ^{(a, a)}}
{\mG ^{(a,\emptyset)}_{aa}}+ \frac{1}{N} \left( \left(   X \wt {\cal S}\cal S X^T  \right)_{aa}+1\right)+ |z|^2}
\\\nonumber
&=\frac{-w(1+(XG^{(a,a)}X^T)_{aa})}{-w(1+m ^{(a, a)})(1+m ^{(a, a)}+Z_{a}^{(a)})
+  \frac1N\left( \left(  X \wt {\cal S}\cal S X^T \right)_{aa}+1\right)+ |z|^2}
\end{align}
We write this denominator as
\be\label{jjssm}
\left(-w(1+m ^{(a, a)})(1+m ^{(a, a)})+|z|^2\right) +\left(-w(1+m ^{(a, a)}) Z_{a}^{(a)} 
+  \frac1N\left( \left(   X \wt {\cal S}\cal S X^T \right)_{aa}+1\right) \right)
 \ee
 With \eqref{230}, \eqref{longsim}, \eqref{53snew}, we can bound the first term in the second bracket as follows:
 $$
\chi_a |w(1+m ^{(a, a)}) Z_{a}^{(a)} |\leq N^{-\e/5}
 $$
 holds with 1-high probability. Together with \eqref{xsx2} and \eqref{abcB2}, with 1-high probability, we can bound the second bracket of \eqref{jjssm} as  
\be\label{jhqd0}
\chi_a \left(-w(1+m ^{(a, a)}) Z_{a}^{(a)} 
+  \frac1N\left( \left(   X \wt {\cal S}\cal S X^T  \right)_{aa}+1\right) \right)
\le N^{- \e/6}
\ee
On the other hand, we claim for some $C>0$, the following inequality holds with 1-hight probability. 
\be\label{jhqd}
\chi_a\left|-w(1+m ^{(a, a)})(1+m ^{(a, a)} ) + |z|^2\right|\ge \chi_a(\log N)^{-C}
\ee
 If \eqref{jhqd} does not hold,  then $\chi_a=1$ and $1+m^{(a,a)}=(-|z|+O(\log N)^{-C})w^{-1/2}$. With \eqref{boundmGm}, \eqref{229} and $||z|-1|\le 2\e$ , we obtain
 \be\label{hs1}
1+  m_\mG^{(a, \emptyset)}=(-|z|+O(\log N)^{-C})w^{-1/2},\ee
It follows from $1+m^{(a,a)}=(-|z|+O(\log N)^{-C})w^{-1/2}$ and  \eqref{231} that 
\be\label{hs2}
\quad  \mG^{(a, \emptyset)} _{aa} =( |z|^{-1}+O(\log N)^{-C})w^{-1/2}
\ee
 Inserting them into \eqref{130}, with \eqref{wzzln},  we have 
\be\label{hs3}
|\cal Z_a|= O(\log N)^{-C}w^{-1/2}
\ee

Now insert \eqref{hs1}, \eqref{hs2} and \eqref{hs3} into \eqref{110},  we obtain  $|G_{aa}|\ge (\log N)^{C-1}|w|^{-1/2}$ for any $C>0$, which contradacts \eqref{53s}. Therefore, \eqref{jhqd} must hold for some $C>0$. 

Recall the denominator of the r.h.s. of \eqref{ydnn} equals to the sum of the l.h.s. of  \eqref{jhqd0},\eqref{jhqd} (see \eqref{jjssm}). Then inserting \eqref{jhqd0},\eqref{jhqd}  into  \eqref{jjssm},  
  we have that for any fixed $D$, there exists $C_{\e, D}$, such that with $1$-high probability,  
\begin{align} \label{tstz}
 \frac{\chi_a}{1+  m_\mG^{(a, \emptyset)}+   |z|^2 \mG^{(a, \emptyset)}_{aa} }&=-\chi_aw(1+(XG^{(a,a)}X^T)_{aa}) 
\\\nonumber
 &*\sum_{k=1}^{C_{\e, D}}\frac{\left( ( -w(1+m ^{(a, a)} )Z^{(a)}_a +\frac1N\left( \left(   X \wt {\cal S}\cal S X^T  \right)_{aa}+1\right)\right)^{k-1}
}{\left(-w(1+m ^{(a, a)})(1+m ^{(a, a)} ) + |z|^2\right)^k}+O_\prec(N^{-D})
 \end{align}
For the  terms in \eqref{tstz}, we apply \eqref{qbq2}  on  $(XG^{(a,a)}X^T)_{aa}$ and $Z^{(a)}_a$,  apply \eqref{232} on $(1+m^{(a,a)})$, apply  \eqref{xsx2} on $\left(   X \wt {\cal S}\cal S X^T  \right)_{aa}$, apply  \eqref{abcB2} on $\gamma$ and apply  \eqref{jhqd} on the denominator of \eqref{tstz}, we obtain 
\begin{align} \nonumber
 \frac{\chi_a}{1+  m_\mG^{(a, \emptyset)}+   |z|^2 \mG^{(a, \emptyset)}_{aa} }&
 \in_n
 -\chi_a
 \left(w\cal F+w\alpha\cal F\right) 
 \sum_{k=1}^{C_{\e, D}}  \left(   -w\alpha^2\cal F +\gamma^2\cal F \right)^{k-1}
 +O_\prec(N^{-D})
 \end{align}
With  the bounds of $\alpha$ and $\gamma$ in  \eqref{abcB}, \eqref{abcB2} and \eqref{abcB4},  it  implies \eqref{DDs3}. Combining \eqref{zGF}, \eqref{DDs3} and \eqref{tjjs}, we obtain \eqref{cGaa}.

\bigskip

{\bf 3.} {For \eqref{cGaab}}, it clearly follows  the Def.  \ref{defSS}, \eqref{boundS} and  Def. \ref{def:Fk}. 


\bigskip

{\bf 4.} {Now we prove \eqref{cGab}}. First with \eqref{110c} and \eqref{1328}, we have 
\be\label{mtmt}
 (\mG_{aa})^{-1}=-w(1+ (YG^{(\emptyset,a)}Y^T)_{aa}).
\ee Applying  $\eqref{Gik}$ on $G_{ak}$ with $i=a$, recalling $Y=X-zI$, we have 
\begin{align}\label{mtmt2}
G_{ak}=&G^{(\emptyset,a)}_{ak}
+w\mG_{aa}\left((G^{(\emptyset,a)}X^T)_{aa}-z^* G^{(\emptyset,a)}_{aa}\right)
\left( (XG^{(\emptyset,a)})_{ak}-z G^{(\emptyset,a)}_{ak}\right) \\\nonumber
 =&G^{(\emptyset,a)}_{ak}
+ w\mG_{aa} G^{(\emptyset,a)}_{aa}  |z|^2 G_{ak}^{(\emptyset,a)}    
-zwG_{ak}^{(\emptyset,a)}  \mG_{aa} (G^{(\emptyset,a)}X^T)_{aa}
\\\nonumber
  &-z^* w\mG_{aa} G^{(\emptyset,a)}_{aa} (XG^{(\emptyset,a)})_{ak}
    +w\mG_{aa} (G^{(\emptyset,a)}X^T)_{aa} (XG^{(\emptyset,a)})_{ak}
  \end{align}
Writing the first term in the   r.h.s. as $G^{(\emptyset,a)}_{ak}\mG_{aa}(\mG_{aa})^{-1}$ and applying  \eqref{mtmt} on $(\mG_{aa})^{-1}$, 
we can write  the first three terms in the  r.h.s. of \eqref{mtmt2} as
 \be\nonumber
   \left(\!- 1\!- (XG^{(\emptyset,a)}X^T)_{aa}+ z^*(XG^{(\emptyset,a)})_{aa}\right)w\mG_{aa}G^{(\emptyset,a)}_{ak}
\ee
Therefore 
\begin{align}\label{qpys}
G_{ak}=&\left(\!- 1\!- (XG^{(\emptyset,a)}X^T)_{aa}+ z^*(XG^{(\emptyset,a)})_{aa}\right)w\mG_{aa}G^{(\emptyset,a)}_{ak}
+\left(-z^*  G^{(\emptyset,a)}_{aa}  
    + (G^{(\emptyset,a)}X^T)_{aa}\right)w\mG_{aa} (XG^{(\emptyset,a)})_{ak}
  \end{align}
  Inserting  \eqref{cGaaa}-\eqref{cGaab}, \eqref{23cl}, \eqref{24cl}, the fact: $\alpha\beta=\chi_a$ and \eqref{xsx3} into \eqref{qpys},  we have 
  \begin{align} \nonumber
\chi_aG_{ak}\in_n&\left(1+\alpha+\beta\gamma+\beta\gamma^2\right)
\cal F^{[\emptyset]} 
\left(\gamma  \cal F_{1/2}^{[k]}
 + \cal F_{0,X}^{[k]} 
\right)
 +  \left(1+ \gamma\right)\cal F^{[\emptyset]}(XG^{(\emptyset,a)})_{ak} +O_\prec(N^{-D})
  \end{align}
More precisely, here what we  used is the $G$-version  of \eqref{23cl}, \eqref{24cl}, i.e., 
$$\chi_a(XG^{(\emptyset,a)})_{aa}\in_n \beta\gamma \cal F \quad {\rm and  }\quad  
\chi_a ( X   G^{( \emptyset,a)}X^T )_{aa}\in( \alpha+\beta\gamma^2)\cal F.
$$
They follows from \eqref{23cl}, \eqref{24cl} and the symmetry between $G$ and $\mG$. 

Next using  \eqref{abcB4}, we have 
   \begin{align}\label{qpys2}
\chi_aG_{ak}\in_n&\left( \alpha +\beta\gamma  \right)
\cal F^{[\emptyset]} 
\left(\gamma  \cal F_{1/2}^{[k]}
 +  \cal F_{0,X}^{[k]}
\right)
 +   \cal F^{[\emptyset]} (XG^{(\emptyset,a)})_{ak} +O_\prec(N^{-D})\\\nonumber
 \in_n& 
\chi_a \sqrt{\frac{\alpha}{N\eta}}\cal F_{1/2}^{[k]}\cdot \cal  F^{[\emptyset]} 
 +
 \left(\alpha +\beta\gamma \right)\cal F^{[\emptyset]}\cdot \cal F_{0,X}^{[k]} 
 + \chi_a(XG^{(\emptyset,a)})_{ak} \cal F^{[\emptyset]} 
 +O_\prec(N^{-D})
  \end{align}
 For $ (XG^{(\emptyset,a)})_{ak}$ in \eqref{qpys2}, using \eqref{111}, for $ k\neq a$   we have (note: $s$ can be $a$)
 $$
 G_{sk}^{(a,a)}= G_{sk}^{( \emptyset,a)}-\frac{ G_{sa}^{( \emptyset,a)} G_{ak}^{( \emptyset,a)}}{ G_{aa}^{( \emptyset,a)}  }
$$
Together with \eqref{xsx3}, \eqref{cGaaa} and \eqref{cGaab}, it implies that 
   \begin{align}\label{jjys}
  \chi_a (XG^{(\emptyset,a)})_{ak}&=\chi_a (XG^{(a,a)})_{ak}
   +\chi_a\frac{(XG^{(\emptyset,a)})_{aa}}{G^{(\emptyset,a)}_{aa}}G^{(\emptyset,a)}_{ak} 
 \\\nonumber
 &\in _n \chi_a(XG^{(a,a)})_{ak}
   + \beta\gamma^2 \cal  F^{[\emptyset]}\cdot \cal F_{1/2}^{[k]} +  \beta\gamma\cal F^{[\emptyset]}\cdot \cal F_{0,X}^{[k]}+O_\prec(N^{-D})
 \end{align}
It follows from \eqref{qbq}, (note: $ |w^{-1/2}|\gamma=\alpha^{1/2} (N\eta) ^{-1/2} $)  that 
 $$\chi_a (XG^{(a,a)})_{ak}\in_n\sqrt{\frac{\alpha}{N\eta}}\cal F_{1/2}^{[k]} +\alpha \cal F_0^{[k]}X_{ak}\in_n
 \sqrt{\frac{\alpha}{N\eta}}\cal F_{1/2}^{[k]} +\alpha \cal F_{0,X}^{[k]}.$$
  Inserting it into \eqref{jjys}, with Lemma \ref{yzb}, 
  we obtain 
   \begin{align}\label{jaw}
 \chi_a  (XG^{(\emptyset,a)})_{ak} 
  \in_n& \sqrt{\frac{\alpha}{N\eta}}\cal F_{1/2}^{[k]}\cdot \cal  F^{[\emptyset]} 
 +
( \alpha +\beta\gamma) \cal F^{[\emptyset]}\cdot \cal F_{0,X}^{[k]} 
  +O_\prec(N^{-D})
  \end{align} 
 Together with  \eqref{qpys2},  we obtain \eqref{cGab}.

  \bigskip
  
 {\bf 5.} {Now we prove \eqref{cGkl}}. With \eqref{jaw}, \eqref{cGaab} and Lem.  \ref{yzb}, we have
 \be\label{busor}
 \chi_\alpha\left((G^{(\emptyset,a)}X^T)_{ka}-z^* G^{(\emptyset,a)}_{ka}\right)\in_n 
 \beta\gamma\cal F^{[k]}_{1/2}\cdot \cal F^{[\emptyset]}+\beta \cal F^{[\emptyset]} \cdot\cal F_{0,X}^{[k]} +O_\prec(N^{-D})
 \ee
 Together with  \eqref{Gik}, \eqref{mtmt}, \eqref{cGaaa} and \eqref{abcB3},   we can  write $G_{kl}$ as follow,
\begin{align}\label{hzz23}
\chi_a\left(G_{kl}-G^{(\emptyset,a)}_{kl}\right)
=& \chi_aw\mG_{aa}\left((G^{(\emptyset,a)}X^T)_{ka}-z^* G^{(\emptyset,a)}_{ka}\right)
\left( (XG^{(\emptyset,a)})_{al}-z G^{(\emptyset,a)}_{al}\right) 
 \\\nonumber
\in_n& \frac1\beta\cal F^{[\emptyset]}
\left(\beta\gamma\cal F^{[k]}_{1/2}\cdot \cal F^{[\emptyset]}
+\beta \cal F^{[\emptyset]}\cdot \cal F_{0,X}^{[k]} 
\right)
\left(\beta\gamma\cal F^{[l]}_{1/2}\cdot \cal F^{[\emptyset]}
+
\beta \cal F^{[\emptyset]} \cdot\cal F_{0,X}^{[l]} 
\right) 
+O_\prec(N^{-D})
\\\nonumber
\in_n &  \frac{1}{N\eta} \cal F^{[k]}_{1/2} \cdot \cal F^{[l]}_{1/2} \cdot\cal F^{[\emptyset]}
+\beta \gamma \left( \cal F_{0,X}^{[k]} \cdot \cal F_{1/2}^{[l]}
+
\cal F_{0,X}^{[l]} \cdot\cal F_{1/2}^{[k]}\right)  \cdot \cal F^{[\emptyset]}  
+\beta \cal F_{0,X}^{[k]}\cdot \cal F_{0,X}^{[l]} \cdot\cal F^{[\emptyset]}
 +O_\prec(N^{-D})\end{align}
Furthermore, with \eqref{111}, \eqref{cGaaa}, \eqref{cGaab} and \eqref{abcB3}, we can write $G^{(\emptyset,a)}_{kl}$ as
\begin{align}\nonumber
 \chi_a\left(G^{(\emptyset,a)}_{kl} -G^{(a,a)}_{kl}\right)
&= \chi_a\frac{G^{(\emptyset,a)}_{ka}G^{(\emptyset,a)}_{al}}{G^{(\emptyset,a)}_{aa}}
=
  \chi_aG^{(\emptyset,a)}_{aa}\frac{G^{(\emptyset,a)}_{ka}}{G^{(\emptyset,a)}_{aa} }\,\frac{G^{(\emptyset,a)}_{al}}{G^{(\emptyset,a)}_{aa} }
\\\nonumber
&\in_n  \frac{\chi_a}{N\eta} \cal F^{[k]}_{1/2}  \cal F^{[l]}_{1/2} \cal F^{[\emptyset]}
+\beta \gamma \left( \cal F_{0,X}^{[k]}  \cal F_{1/2}^{[l]}
+
\cal F_{0,X}^{[l]} \cal F_{1/2}^{[k]}\right)   \cal F^{[\emptyset]}  
+\beta \cal F_{0,X}^{[k]} \cal F_{0,X}^{[l]} \cal F^{[\emptyset]}
+O_\prec(N^{-D})
 \end{align}
 Therefore, together with \eqref{hzz23}, we obtain  \eqref{cGkl}. 

\qed

Next, we write the terms appeared in   the Lemma \ref{AYFF}  as polynomials in $\cal F$, $\cal F_{1/2}$ and $\cal F_{1/2}\cdot \cal F$ (with proper coefficients   and ignorable error terms).
 
\begin{lemma}\label{cbound2}
Let $w\in I_\e$ and $||z|-1|\le2\e$. Under the assumption of Lemma \ref{mainnewL},  for any  fixed large $D>0$, with  $\chi_a$ defined  in \eqref{defchaa} and $F_{ 0, X}^{[k]} $ defined in \eqref{gniz}, we have that for $k\neq a$
 \begin{align} \label{cm}
\chi_a  (m - m^{(a, a)})&\in_n\frac{1}{N\eta}\cal F+O_\prec(N^{-D})
\\\label{cGbb}
\chi_aG_{bb} &\in_n \beta\cal F +O_\prec(N^{-D})
\\\label{cYGaa}
\chi_a(YG )_{aa} 
& \in_n   \cal F +O_\prec(N^{-D})
\\ \label{cYGab}
  \chi_a(YG )_{ak}
  & \in_n  \gamma\cal F^{[k]}_{1/2}\cdot \cal F^{[\emptyset]}+   \cal F_{0,X}^{[k]}    \cal F^{[\emptyset]} + O_\prec(N^{-D}), 
\\\label{cYG2ab}
 \chi_a(YG^2 )_{ab}
 &\in_n \frac{\chi_a}{ \eta}  \cal F_{1/2}\cdot \cal F +X_{ba}\cal F+O_\prec(N^{-D}).  
\\\label{cG2aa}
\chi_a(G^2)_{aa}& \in_n  \frac{ \alpha }{ \eta} \cal F+O_\prec(N^{-D})
\\\label{cG2bb}
\chi_a(G^2)_{bb}&\in_n    \beta  \eta^{-1}\cal F
+O_\prec(N^{-D})
\\\label{cYG2Y}
\chi_a (YG^2Y^*)_{aa}
&\in_n\frac{w\alpha}{\eta}\cal F+O_\prec(N^{-D})
\end{align}
 \end{lemma}

\bigskip

{\it Proof of Lemma \ref{cbound2}: }  {\bf 1.} {For \eqref{cm}}, using \eqref{cGkl} and  \eqref{cGaa}, we have  
\begin{align}\nonumber
\chi_a(m - m^{(a, a)})&=\frac{\chi_a}NG_{aa}
+\chi_a\frac1N\sum_{k\neq a}\left(G_{kk}-G^{(a,a)}_{kk}\right)\\\nonumber
&\in_n 
\frac{\alpha}{N}\cal F
+
 \frac1N\sum_{k\neq a} \Bigg(\frac{\chi_a}{N\eta} \cal F^{[k]}_{1/2}  \cal F^{[k]}_{1/2} 
+\beta \gamma  \cal F_{0,X}^{[k]}  \cal F_{1/2}^{[k]}   +\beta
  \cal F_{0,X}^{[k]}    \cal F_{0,X}^{[k]}  \Bigg)\cal F^{[\emptyset]}+O_\prec(N^{-D})
\\\nonumber
&\in_n \left(\frac{\alpha}{N}+ \frac{1}{N\eta}+\frac{\beta\gamma}N+\frac{\beta }N\right)\cal F
 +O_\prec(N^{-D})
\end{align}
Here for the last $\in _n$, we used 
\be\label{XFFX}
 \sum_{k\neq a}  \frac1N\cal F^{[k]}_{1/2}  \cal F^{[k]}_{1/2} \in \cal F,\quad 
 \sum_{k\neq a}   \cal F_{0,X}^{[k]}    \cal F^{[k]}_{1/2} \in \cal F, \quad  
    \sum_{k\neq a} \cal F_{0,X}^{[k]}    \cal F_{0,X}^{[k]}   \in \cal F.
   \ee  Then with \eqref{abcB} and \eqref{abcB2}, we obtain \eqref{cm}. 
 
 \bigskip
 
  {\bf 2.} {For \eqref{cGbb}}, it follows from \eqref{cGkl}, $\cal F_{1/2}\cdot \cal F_{1/2}\subset \cal F$ and the fact: $X_{ab}=0$ that 
  \begin{align}\nonumber
 \chi _a G_{bb} &\in _n  \chi_aG^{(a,a)}_{bb} +\frac{\chi_a}{N\eta}\cal F+\beta\gamma
 X_{ba}\cal F_{1/2} \cdot\cal F
+\beta X_{ba}X_{ba} \cal F\\\nonumber
&\in_n \alpha \cal F+\frac{\chi_a}{N\eta}\cal F +(\beta+\gamma\beta) \cal F
\end{align}
where we used \eqref{234} on $G^{(a,a)}_{bb}$, $X_{ba}\in \cal F_{1/2}$. Now using Lemma \ref{yzb}, we obtain \eqref{cGbb}. 

{\bf 3.}  {For \eqref{cYGaa}},  with \eqref{Gik} and \eqref{mtmt}, we can write it as \begin{align}\label{erz}
\chi_a(YG )_{aa}=-\chi_aw\mG_{aa}(YG^{(\emptyset, a)} )_{aa}
&=-\chi_aw\mG_{aa} \left((XG^{(\emptyset, a)} )_{aa} - z G_{aa}^{(\emptyset, a)}\right)
 \end{align}
Then with   \eqref{cGaa}, \eqref{cGaaa} and \eqref{23cl}, we obtain \eqref{cYGaa}. 

{\bf 4.} {Now we prove \eqref{cYGab}}, with \eqref{Gik} and \eqref{mtmt} again, we write it is  \begin{align}\label{erz2}
(YG )_{ak}=-w\mG_{aa}(YG^{(\emptyset, a)} )_{ak}
&=-w\mG_{aa} \left((XG^{(\emptyset, a)} )_{ak} - z G_{ak}^{(\emptyset, a)}\right)
 \end{align}
Then using   \eqref{busor} and \eqref{cGaa},    we obtain \eqref{cYGab}. 

{\bf 5.} {For \eqref{cYG2ab}}, by definition, we write $(YG^2 )_{ab}$ as
 \be\nonumber
(YG^2 )_{ab}=\sum_ {k\neq a } (YG)_{ak}G_{kb} +(YG)_{aa}G_{ab} 
\ee
Then using \eqref{cYGab}, \eqref{cGkl},  with $X_{ab}=0$, we get 
\begin{align}\label{5119}
 &\chi_a\sum_ {k\neq a } (YG)_{ak}G_{kb}
 \\\nonumber
 & \in_n   
  \sum_{k\neq a } \cal F^{[\emptyset]}\left( \gamma\cal F^{[k]}_{1/2} +  \cal F_{0,X}^{[k]} \right)\Bigg( G^{(a,a)}_{kb} 
 + \frac{\chi_a}{N\eta} \cal F^{[k]}_{1/2}  \cal F^{[b]}_{1/2} 
+\beta \gamma \left( \cal F_{0,X}^{[k]} \cal F_{1/2}^{[b]} +\cal F_{0,X}^{[b]}\cal F_{1/2}^{[k]}\right) +\beta  \cal F_{0,X}^{[k]} \cal F_{0,X}^{[b]} \Bigg)\end{align}
With \eqref{XFFX} and Lemma \ref{yzb}, we obtain 
\be\nonumber
  \chi_a\sum_ {k\neq a } (YG)_{ak}G_{kb}\in_n   
 \cal F^{[\emptyset]} 
 \sum_{k\neq a } \left( \gamma\cal F^{[k]}_{1/2} + \cal F_{0,X}^{[k]} \right)  G^{(a,a)}_{kb}
 +\frac{\chi_a}{\eta} \left(\cal F^{[b]}_{1/2} \cal F^{[\emptyset]}+\cal F^{[b]} X_{ba}\right) 
\ee
Then applying \eqref{qbq} on $ G^{(a,a)}_{kb}$, we obtain $ \chi_aG^{(a,a)}_{kb}\in \left(|w|^{-1/2}\gamma+\delta_{kb}\alpha\right)\cal F_0^{[k,b]}$. Now with
$$
\sum_{k\neq a}\cal F^{[k ]}_{1/2}\cdot\cal F_0^{[k,b]}\in N \cal F_{1/2},\quad \quad  
\sum_{k\neq a}\cal F^{[k ]}_{0, X}\cdot\cal F_0^{[k,b]}\in   \cal F_{1/2}
$$
and  Lemma \ref{yzb} again,  we get 
\be\nonumber
\chi_a\sum_ {k\neq a } (YG)_{ak}G_{kb}
\in_n
\frac{\chi_a}{\eta} \left(\cal F^{[b]}_{1/2} \cal F^{[\emptyset]}+\cal F^{[b]} X_{ba}\right) 
\ee
Similarly, with \eqref{cYGaa},  \eqref{cGab} and  Lemma \ref{yzb} again, we obtain 
\be\nonumber
\chi_a(YG)_{aa}G_{ab}  \in_n\frac{\chi_a}{\eta} \left(\cal F^{[b]}_{1/2} \cal F^{[\emptyset]}+\cal F^{[b]} X_{ba}\right) 
\ee
and we obtain \eqref{cYG2ab}.

\bigskip

{\bf 6.} {For  \eqref{cG2aa}}, we write $(G^2)_{aa}$ as  
\be\nonumber
\chi_a(G^2)_{aa}=\chi_a\sum_{k\neq a } G_{ak}G_{ka}+ \chi_a(G_{aa})^2\in_n\sum_{k\neq a } \chi_aG_{ak}G_{ka}+  \frac{\alpha}{\eta}\cal F
\ee
where for the second $\in_n$, we used \eqref{cGaa} and \eqref{abcB}.  As in \eqref{5119}, using \eqref{cGab} and \eqref{XFFX}, we have 
\begin{align}\nonumber
\chi_a\sum_{k\neq a } G_{ak}G_{ka}&=\sum_{k\neq a } \Bigg( \sqrt{\frac{\alpha}{N\eta}}  \cal F_{1/2}^{[k]}\cdot \cal F ^{[\emptyset]}
+\left(\alpha + \beta\gamma\right) \cal F^{[\emptyset]}\cdot \cal F_{0,X}^{[k]}\Bigg)^2
\\\nonumber
&\in_n \left(\frac{\alpha}{\eta}+  \left(\alpha + \beta\gamma\right)^2  \right)\cal F\end{align}
Then with Lemma \ref{yzb}, we obtain \eqref{cG2aa}.

{\bf 7.} { \eqref{cG2bb}}, we write it as
\be\nonumber
\chi_a(G^2)_{bb}=\sum_{k\neq a,b } G_{bk}G_{kb}+ (G_{bb})^2+ (G_{ab})^2
\ee
With \eqref{cGkl},  \eqref{XFFX} and $\sum_k \cal F^{[k,b]}_\alpha\in N\cal F^{[b]}_\alpha$, $(\alpha= 0,1/2,\emptyset)$, after a tedious calculation, we get
\begin{align}\nonumber
&\chi_a\sum_{k\neq a,b } G_{bk}G_{kb}\\\nonumber
&=\sum_{k\neq a,b }\left(\sqrt{\frac{\alpha}{N\eta}}\cal F_0^{[k,b]}+\frac{\chi_a}{N\eta} \cal F^{[k]}_{1/2}  \cal F^{[b]}_{1/2} 
+\beta \gamma \left( \cal F_{0,X}^{[k]} \cal F_{1/2}^{[b]}+ \cal F_{0,X}^{[b]} \cal F_{1/2}^{[k]}\right)      +\beta \cal F_{0,X}^{[k]} \cal F_{0,X}^{[b]}\right)^2\cal F^{[\emptyset]}
\\\nonumber
&\in_n 
\frac{\alpha}{\eta}\cal F
+ \left(\frac1\eta\sqrt{\frac{ \alpha }{N\eta}} +\beta\gamma \sqrt{\frac{ \alpha }{N\eta}} 
+\frac{1}{N\eta^2} +\frac{\beta\gamma}{N\eta}+\beta^2\gamma^2\right)
\cal F^{[b]}_{1/2}\cal F^{[b]}_{1/2}\cal F\;+ O_\prec(N^{-D})
\\\nonumber
&\quad
+\left(N\beta \gamma\sqrt{\frac{ \alpha }{N\eta}}+\beta\sqrt{\frac{ \alpha }{N\eta}}+\frac{\beta\gamma}{\eta}+ \frac{\beta}{\eta}+\beta^2\gamma^2+\beta^2\gamma\right)
X_{ba}\cal F^{[b]}_{1/2}\cal F 
+\left(\beta^2\gamma^2N+\beta^2\gamma+\beta^2\right)X_{ba}^2\cal F
 \end{align}
Then using Lemma \ref{yzb}, $\cal F^{[b]}_{1/2}\cal F^{[b]}_{1/2}\in \cal F$, 
$X_{ba}\cal F^{[b]}_{1/2}\in \cal F$ and $X_{ba}^2\in \cal F$, we obtain 
 \be\nonumber
\chi_a\sum_{k\neq a,b } G_{bk}G_{kb} 
\in_n \frac{\beta}{\eta}\cal F+O_\prec(N^{-D})
\ee
 Similarly, using \eqref{cGkl},  and Lemma \ref{yzb} we have
 $$
 \chi_aG_{bb}G_{bb} \in_n \left(\alpha+\frac{1}{N\eta} +\beta\gamma+\beta\right)^2 \cal F+O_\prec(N^{-D})
 \in_n\frac{\beta}{\eta}\cal F+O_\prec(N^{-D})
 $$
Using   \eqref{cGaa},  and Lemma \ref{yzb}  we have
 \be\nonumber
 \chi_aG_{ba}G_{ab} 
\in_n \frac{\beta}{\eta}\cal F+O_\prec(N^{-D})
\ee
which completes the proof of  \eqref{cG2bb}. 
\bigskip

{\bf 8.} For  {\eqref{cYG2Y}}, it follows from 
\be\nonumber
 (YG^2Y^T)_{aa} = \mG_{aa}  +w(\mG^2)_{aa} 
\ee
and \eqref{cGaa} and \eqref{cG2aa}. 

\qed 

Now we are ready to prove Lemma \ref{AYFF}, which is  the key lemma in the proof of  our main result. 
 
 \subsection{Proof of lemma \ref{AYFF}. }
First with $m-m^{(a,a)}=O(N\eta)^{-1}$ (see \eqref{boundmGm}) and the definition of $\chi_a$, for any fixed $D>0$, with 1-high probability, we can write the $h(t_X)$ as
\be\nonumber
h(t_X)=\chi_a h(t_X)=\sum_{k=0}^{C_{\e, D}} \frac{1}{k!}   h^{(k)}(t_{X^{(a,a)}})  \chi_a\left(\frac{\re m-\re m^{(a,a)}}{N^{\e} (N\eta)^{-1}}\right)^k+O(N^{-D})
\ee
where constant $C_{\e, D}$   depends on $\e$ on $D$,  and $h^{(k)}$ is the $k-th$ derivative of $h$. Using \eqref{cm} and the fact that $h$ is smooth and supported in $[1,2]$, we obtain 
\be\label{280a}
h(t_X)\in _n \cal F+O_\prec(N^{-D})
\ee
and 
\be\label{280aa}
h(t_X)-h(t_{X^{(a,a)}})\in _n N^{-\e}{\bf 1}(|t_{X^{(a,a)}}|\le 2)\cal F+ O_\prec(N^{-D})
\ee  
Note: ${\bf 1}(|t_{X^{(a,a)}}|\le 2)={\bf 1}(|\re m^{(a,a)}|\le 2N^{\e}(N\eta)^{-1})$. Similarly, one can prove 
\be\label{280b}
h'(t_X), \quad h''(t_X), \quad h'''(t_X) \in _n {\bf 1}(|t_{X^{(a,a)}}|\le 2)\cal F+O_\prec(N^{-D})
\ee
Using \eqref{280a}, \eqref{280aa} and \eqref{cm}, we have 
\begin{align}\nonumber
\left(h(t_X)\re m-h(t_{X^{(a,a)}}) \re  m^{(a,a)}\right)&\in_n  \left( h(t_X)  \re m^{(a,a)}-h(t_{X^{(a,a)}}) \re  m^{(a,a)}\right)+ \frac{1}{N\eta }\cal F+O(N^{-D})  \\\label{ser}
 & \in_n \frac{1}{N\eta }\cal F +O_\prec(N^{-D})
\end{align}
It implies \eqref{ygfld}. 

For \eqref{ygfld2}, recall $B_m(X)$ is defined as 
\be\nonumber
B_m(X):= \frac1{m!}(N^{1-\e}\eta )^{ (m-1)}\left(mh^{( m-1)}(t_{X})+h^{( m)}(t_{X})t_{X}\right)\ee
Then using \eqref{280a}, \eqref{280b} and \eqref{cm}, we obtain \eqref{ygfld2}.

  \bigskip
 
Similarly, for  \eqref{ygfld3}, the terms appearing in the definition \eqref{PPP} have been all bounded in   \eqref{cYG2ab}, \eqref{cGaa}, \eqref{cG2bb}, \eqref{cYG2ab}, \eqref{cYG2Y} and \eqref{cGbb}. With a simple calculation, one can obtain \eqref{ygfld3} and complete the proof.  
      
      \qed

 
 

\end{document}